\documentclass[11pt]{amsart}
\usepackage{amsmath,amsthm,amscd,amssymb}
\usepackage{pb-diagram}
\usepackage{graphicx}
\usepackage{rotating}
\usepackage{xcolor}
\usepackage{amssymb}
\newcommand{\Z}{\mathbb{Z}}
\newcommand{\sP}{\scriptscriptstyle{\mathcal{P}}}
\newcommand{\K}{\mathcal{K}}

\newtheorem{thm}{Theorem}[section]
\newtheorem{lem}[thm]{Lemma}
\newtheorem{prop}[thm]{Proposition}
\newtheorem{cor}[thm]{Corollary}

\theoremstyle{definition}
\newtheorem{ex}[thm]{Example}
\newtheorem{defn}[thm]{Definition}

\numberwithin{equation}{section}
\numberwithin{figure}{section}

\newcommand{\A}{\mathcal{A}}
\newcommand{\G}{\Gamma}

\newcommand{\Q}{\mathbb{Q}}
\newcommand{\R}{\mathcal{R}}

\usepackage{pstricks}
\usepackage{pst-node}

\newcommand{\sss}{\scriptscriptstyle}

\newcommand{\gn}{\ensuremath{G^{\sss (n)}}}
\newcommand{\gnp}{\ensuremath{G^{\sss (n+1)}}}

\title[$2$-Torsion in the $n$-Solvable Filtration of the Knot concordance group]{2-Torsion in the $\boldsymbol{n}$-Solvable Filtration of the Knot concordance group}

\author{Tim D. Cochran$^{\dag}$}
\address{Department of Mathematics MS-136, P.O. Box 1892, Rice University, Houston, TX 77251-1892}
\email{cochran@rice.edu}

\author{Shelly Harvey$^{\dag\dag}$}
\address{Department of Mathematics MS-136, P.O. Box 1892, Rice University, Houston, TX 77251-1892}
\email{shelly@rice.edu}

\author{Constance Leidy$^{\dag\dag\dag}$}
\address{Department of Mathematics, Wesleyan University, Wesleyan Station, Middletown, CT 06459}
\email{cleidy@wesleyan.edu}

\thanks{\noindent $^{\dag}$Partially supported by NSF DMS-0706929}
\thanks{ $^{\dag\dag}$Partially supported
by NSF CAREER DMS-0748458}
\thanks{$^{\dag\dag\dag}$Partially supported by NSF DMS-0805867}


\begin{document}

\begin{abstract} Cochran-Orr-Teichner introduced in ~\cite{COT} a natural filtration of the smooth knot concordance group $\mathcal{C}$
$$
\cdots \subset \mathcal{F}_{n+1} \subset \mathcal{F}_{n.5}\subset\mathcal{F}_{n}\subset\cdots \subset
\mathcal{F}_1\subset \mathcal{F}_{0.5} \subset \mathcal{F}_{0} \subset \mathcal{C},
$$
called the ($n$)-solvable filtration. We show that each associated graded abelian group $\{\mathbb{G}_n=\mathcal{F}_{n}/\mathcal{F}_{n.5}~|~n\in\mathbb{N}\}$, $n\geq 2$ contains infinite linearly independent sets of elements of order $2$ (this was known previously for $n=0,1$). Each of the representative knots is negative amphichiral, with vanishing $s$-invariant, $\tau$-invariant, $\delta$-invariants and Casson-Gordon invariants. Moreover each is slice in a rational homology $4$-ball.
In fact we show that there are many distinct such classes in $\mathbb{G}_n$, distinguished by their Alexander polynomials and, more generally, by the torsion in their higher-order Alexander modules.
\end{abstract}

\maketitle
\section{Introduction}\label{intro}

A (classical) knot $K$ is the image of a smooth embedding of an oriented circle in $S^3$. Two knots, $K_0\hookrightarrow S^3\times \{0\}$ and $K_1\hookrightarrow S^3\times \{1\}$, are concordant if there exists a proper smooth embedding of an annulus into $S^3\times [0,1]$ that restricts to the knots on $S^3\times \{0,1\}$. Let $\mathcal{C}$ denote the set of (smooth) concordance classes of knots. The equivalence relation of concordance first arose in the early $1960$'s in work of Fox, Kervaire and Milnor in their study of isolated singularities of $2$-spheres in $4$-manifolds and, indeed, certain concordance problems are known to be \emph{equivalent} to whether higher-dimensional surgery techniques ``work'' for topological $4$-manifolds \cite{FM}\cite{KM1}\cite{CF}. It is well-known that $\mathcal{C}$ can be endowed  with the structure of an abelian group (under the operation of connected-sum), called the smooth knot concordance group. The identity element is the class of the trivial knot. Any knot in this class is concordant to a trivial knot and is called a slice knot. Equivalently, a slice knot is one that is the boundary of a smooth embedding of a $2$-disk in $B^4$. In general, the abelian group structure of $\mathcal{C}$ is still poorly understood.  But much work has been done on the subject of knot concordance (for excellent surveys see ~\cite{Go1} and ~\cite{Li1}). In particular,  ~\cite{COT} introduced a natural filtration of $\mathcal{C}$ by subgroups
$$
\cdots \subset \mathcal{F}_{n+1} \subset \mathcal{F}_{n.5}\subset\mathcal{F}_{n}\subset\cdots \subset
\mathcal{F}_1\subset \mathcal{F}_{0.5} \subset \mathcal{F}_{0} \subset \mathcal{C}.
$$
called the ($n$)-\emph{solvable filtration} of $\mathcal{C}$ and denoted $\{\mathcal{F}_{n}\}$ (defined in Section~\ref{sec:series}).
The non-triviality of $\mathcal{C}$ can be measured in terms of the associated graded abelian groups $\{\mathbb{G}_n=\mathcal{F}_{n}/\mathcal{F}_{n.5}~|~n\in\mathbb{N}\}$ (here we ignore the other ``half'' of the filtration, $\mathcal{F}_{n.5}/\mathcal{F}_{n+1}$, where almost nothing is known). This paper is concerned with elements of order two in $\mathcal{C}$ and, more generally, with elements of order two in $\mathbb{G}_n$.

We will review some of the history of $2$-torsion phenomena in $\mathcal{C}$ in the context of the $n$-solvable filtration. One of the earliest results concerning $\mathcal{C}$ was an epimorphism constructed by Fox and Milnor ~\cite{FM}
$$
FM:\mathcal{C}\twoheadrightarrow \Z_2^\infty.
$$
Soon thereafter, Levine constructed an epimorphism
\begin{equation}\label{eq:algconcgrp}
\mathcal{C}\twoheadrightarrow \mathcal{AC}\cong \Z^\infty\oplus\Z_2^\infty\oplus \Z_4^\infty,
\end{equation}
to a group, $\mathcal{AC}$, that became known as the \emph{algebraic knot concordance group}. Any knot in the kernel of \eqref{eq:algconcgrp} is called an \textbf{algebraically slice} knot. In terms of the $n$-solvable filtration, Levine's result is \cite[Remark 1.3.2, Thm. 1.1]{COT}:
$$
\mathbb{G}_0\cong \Z^\infty\oplus\Z_2^\infty\oplus \Z_4^\infty.
$$
It is known that there exist elements of order two in $\mathcal{C}$ that realize some of the above $2$-torsion invariants. Let $\overline{K}$ denote the \textbf{mirror image} of the oriented knot $K$, obtained as the image of $K$ under an orientation reversing homeomorphism of $S^3$; and let $r(K)$ denote the \textbf{reverse} of $K$, which is obtained by merely changing the orientation of the circle. Then it is known that $K\# r(\overline{K})$ is a slice knot, so the inverse of $[K]$ in $\mathcal{C}$, denoted $-[K]$, is represented by $r(\overline{K})$, denoted $-K$. A knot $K$ is called \textbf{negative amphichiral} if $K$ is isotopic to $r(\overline{K})$. It follows that, for any negative amphichiral knot $K$, $K\# K$ is a slice knot, since it is isotopic to $K\# -K$. Hence negative amphichiral knots represent elements of order either $1$ or $2$ in $\mathcal{C}$. It is a conjecture of Gordon that every class of order two in $\mathcal{C}$ can be represented by a negative amphichiral knot \cite{Go1}.

In fact the work of Milnor and Levine in the $1960$'s resulted in a more precise statement:
$$
\mathbb{G}_{0}\cong\bigoplus_{\substack{p(t)}}\left (\Z^{r_p}\oplus\Z_2^{m_p}\oplus \Z_4^{n_p}\right )
$$
where the sum is over all primes $p(t)\in \Z[t]$ where $p(t)\doteq p(t^{-1})$ and $p(1)=\pm 1$ \cite[Sections 10,11,24]{Le10}\cite{Sto}\cite[p.131]{Hi}. That is, the algebraic concordance group (and $\mathbb{G}_{0}$) admits a certain \emph{$p(t)$-primary decomposition}, wherein a knot has a nontrivial $p(t)$-primary part only if $p(t)$ is a factor of its Alexander polynomial. (Indeed, Levine and Stoltzfus classified $\mathbb{G}_{0}$ by first splitting the Witt class of the Alexander module (with its Blanchfield form) into its $p(t)$-primary parts).

In the $1970$'s the introduction of Casson-Gordon invariants in ~\cite{CG1}\cite{CG2} led to the discovery that the subgroup of algebraically slice knots was of infinite rank and contained infinite linearly independent sets of elements of order two \cite{Ji1}\cite{Li6}. In terms of the $n$-solvable filtration this implies the existence of
$$
\Z^\infty\oplus \Z_2^\infty\subset \mathbb{G}_{1}.
$$
Different $\Z^\infty$-summands were exhibited in \cite{Ki1}\cite{Fr2}. More recent work of Se-Goo Kim ~\cite{KiS} on the ``polynomial splitting'' properties of Casson-Gordon invariants led to a generalization analogous to the result of Milnor-Levine:
$$
\bigoplus_{\substack{p(t)}}\Z^\infty\subset\mathbb{G}_{1}.
$$
Thus there is evidence that $\mathbb{G}_{1}$ also exhibits a $p(t)$-primary decomposition. Further strong evidence is given in \cite{KiKi}. Although a similar statement for the $2$-torsion in $\mathbb{G}_{1}$ has not appeared, it is expected from combining the work of ~\cite{KiS} and Livingston \cite{Li6}. Several authors have shown that certain knots that projected to classes of order $2$ and $4$ in $\mathcal{AC}$ are in fact of \emph{infinite order} in $\mathcal{C}$ \cite{LiN1}\cite{LiN2}\cite{JaN}\cite{GRS}\cite{Lis1}. A number of papers have addressed the non-triviality of $\{\mathbb{G}_{n}\}$, \cite{Gi3}\cite{GL2}\cite{Ki1}\cite{Fr2}\cite{COT}\cite{COT2}\cite{CT},
culminating in ~\cite{CHL3} where it was shown that, for any integer $n$, there exists
$$
\Z^\infty\subset \mathbb{G}_{n}.
$$
Moreover the recent work \cite{CHL5} of the authors resulted in a generalization of the latter fact, along the lines of the Levine-Milnor primary decomposition and ~\cite{KiKi}: for each ``distinct'' $n$-tuple $\mathcal{P}=(p_1(t),...,p_n(t))$ of prime polynomials with $p_{i}(1)=\pm1$, there is a distinct subgroup $\Z^\infty\cong \mathcal{I}(\mathcal{P})\subset\mathbb{G}_n$, yielding a subgroup
\begin{equation}\label{eq:bigsubgroup}
\bigoplus_{\substack{\mathbb{P}_n}}\Z^\infty\cong\bigoplus_{\substack{\mathcal{P}\in\mathbb{P}_n}}\mathcal{I}(\mathcal{P})\subset \mathbb{G}_n.
\end{equation}
Given a knot $K$, such an $n$-tuple encodes the orders of certain submodules of the sequence of higher-order Alexander modules of $K$. Thus one can distinguish concordance classes of knots not only by their classical Alexander polynomials, but also, loosely speaking, \emph{by their higher-order Alexander polynomials}. This result indicates that $\mathbb{G}_{n}$ decomposes not just according to the prime factors of the classical Alexander polynomial, but also according to types of torsion in the higher-order Alexander polynomials.

Here we show corresponding results for $2$-torsion. That is, for any $n\geq 2$, not only will we exhibit
\begin{equation}\label{eq:littletorsionsubgroup}
\Z_2^\infty\subset \mathbb{G}_{n},
\end{equation}
but we also will exhibit many distinct such subgroups
\begin{equation}\label{eq:torsionbigsubgroup}
\bigoplus_{\substack{\mathbb{P}_{n-1}}}\Z_2^\infty\subset \mathbb{G}_n,
\end{equation}
parametrized by their Alexander polynomials and the types of torsion in the higher-order Alexander polynomials. The representative knots are distinguished by families of von Neumann signature defects associated to their classical Alexander polynomials and ``higher-order Alexander polynomials''. The precise statement is given in Theorem~\ref{thm:bigsum}. Each of these concordance classes has a negative amphichiral representative that is smoothly slice in a rational homology $4$-ball. Thus the classical signatures and the Casson-Gordon signature-defect obstructions \cite{CG1} (indeed all metabelian obstructions) vanish for these knots \cite[Theorem 9.11]{COT}. In addition, the $s$-invariant of Rasmussen \cite{Ra1}, the $\tau$-invariant of Ozsv\'{a}th-Szab\'{o} \cite{OS1}, and the $\delta_{p^n}$ invariants of Manolescu-Owens and Jabuka \cite{MO1}\cite{Ja}\cite{OS2} vanish on these concordance classes, since each of these invariants induces a \emph{homomorphism} $\mathcal{C}\to\Z$ and so must have value zero on classes representing torsion in $\mathcal{C}$. Our examples are inspired by those of Livingston, who provided examples that can be used to establish \eqref{eq:littletorsionsubgroup} in the case $n=1$ ~\cite{Li6}. His examples are distinguished by their Casson-Gordon signature defects. Our examples are distinguished by higher-order $L^{(2)}$-signature defects. It is striking that elements of finite order can sometimes be detected by signatures. The key observation is that, unlike invariants such as the classical knot signatures, the $s$ invariant, the $\tau$-invariant, or the $\delta$-invariants, the invariants arising from higher-order signature defects (including Casson-Gordon invariants) are not additive under connected sum. Therefore there is no reason to expect that they would vanish on elements of finite order.

Our work is further evidence that $\mathbb{G}_{n}$ exhibits some sort of primary decomposition, but wherein not only the classical Alexander polynomial, but also some higher-order Alexander polynomials are involved.

We remark that \cite{COT} also defined a filtration, $\{\mathcal{F}_n^{top}\}$, of the topological concordance group, $\mathcal{C}^{top}$. Since it is known, by work of Freedman and Quinn, that a knot lies in $\mathcal{F}_n^{top}$ if and only if it lies in $\mathcal{F}_n$, all of the results of this paper apply equally well, without change, to the filtration $\{\mathcal{F}_n^{top}\}$. Therefore, for simplicity, in this paper we will always work in the smooth category.

\section{The examples}\label{sec:examples}

Our examples are inspired by those of Livingston ~\cite{Li6}, who exhibited an infinite ``linearly independent'' set of negative amphichiral algebraically slice knots. His examples can be used to establish the existence of the aforementioned
$$
\Z_2^\infty\subset \mathbb{G}_{1}.
$$

\begin{subsection}{The Building Blocks}\label{subsec:building blocks}

Consider the knot shown on the left-hand side of Figure~\ref{fig:generalamphi}.
\begin{figure}[htbp]
\setlength{\unitlength}{1pt}
\begin{picture}(367,132)
\put(0,0){\includegraphics{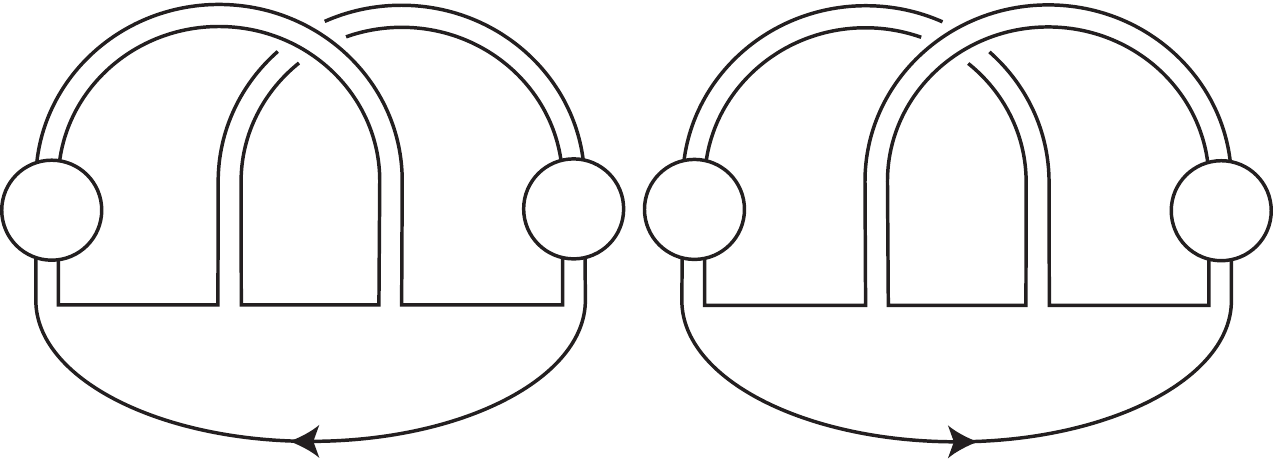}}
\put(155,69){$f(\overline{J})$}
\put(10,69){$J$}
\put(342,69){$f(J)$}
\put(196,69){$\overline{J}$}
\end{picture}
\caption{Families of Negative Amphichiral Knots $K$}\label{fig:generalamphi}
\end{figure}
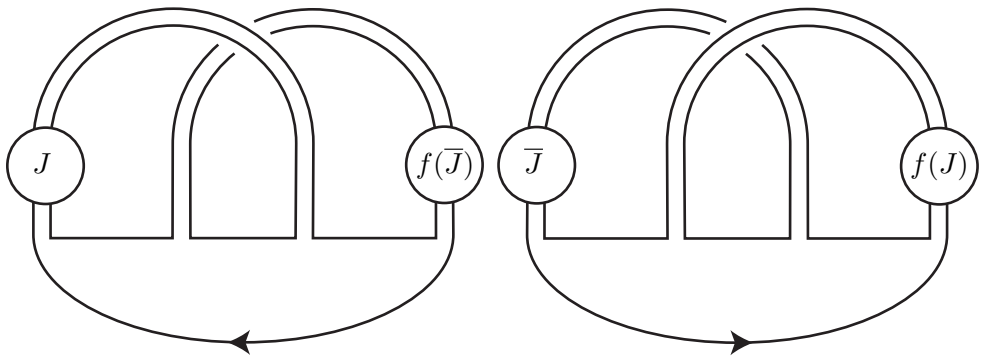
Here $J$ is an arbitrary pure two component string link \cite{LD}\cite{HL1}. The disk containing the letter $J$ symbolizes replacing the trivial $2$-string link by the $2$-string link $J$. Viewing the knot diagram as being in the $xy$-plane ($y$ being vertical), the mirror image can be defined as the image under the reflection $(x,y,z)\mapsto (x,y,-z)$, which alters a knot diagram by replacing all positive crossings by negative crossings and vice-versa. Recall that the image of $J$ under this reflection is denoted $\overline{J}$. We also consider a ``flip'' homeomorphism of $S^3$ which flips over a diagram, given by rotation of $180$ degrees about the $y$-axis or $f(x,y,z)=(-x,y,-z)$. Note that these homeomorphisms commute. Special cases of the following elementary observation appeared in ~\cite[Lemma 2.1]{Li6} ~\cite[p.326]{Li1} and ~\cite[p. 60]{Cha2}.

\begin{lem}\label{lem:amphi} Suppose $J$ is an arbitrary pure two component string link. Then the knot $K$ on the left-hand side of Figure~\ref{fig:generalamphi} is negative amphichiral.
\end{lem}

\begin{figure}[htbp]
\setlength{\unitlength}{1pt}
\begin{picture}(376,168)
\put(0,0){\includegraphics{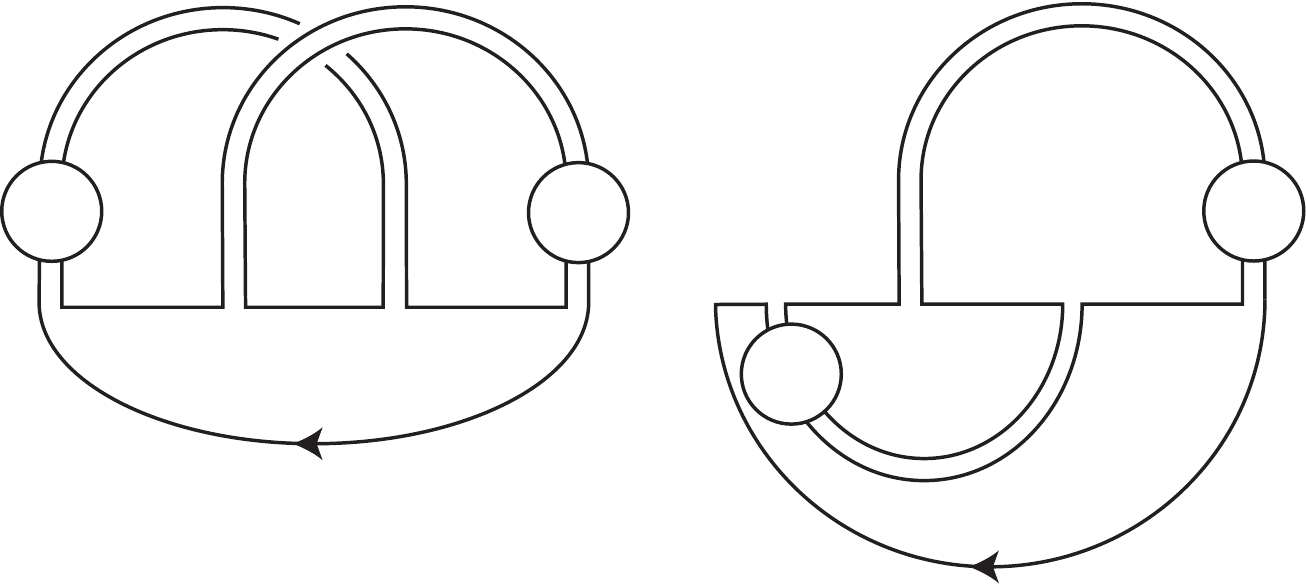}}
\put(352,104){$f(\overline{J})$}
\put(9,104){$J$}
\put(156,104){$f(\overline{J})$}
\put(224,65){\rotatebox{185}{\reflectbox{$J$}}}
\end{picture}
\caption{}\label{fig:amphiproof}
\end{figure}
\begin{proof} The knot on the right-hand side of Figure~\ref{fig:generalamphi} is a diagram for $r(\overline{K})$, since it is obtained by a reflection, in the plane of the paper, of the diagram for $K$, followed by a reversal of the string orientation. Here we use that $f$ commutes with the reflection. We claim that the result is isotopic to $K$. Flipping the diagram (rotating by $180$ degrees about the vertical axis in the plane of the paper), we arrive at the diagram shown on the left-hand side of Figure~\ref{fig:amphiproof}. This is identical to the original diagram of $K$ except that the left-hand band passes under the right-hand band instead of over. But the left-hand band can be ``swung'' around by an isotopy as suggested in the right-hand side of Figure~\ref{fig:amphiproof}, bringing it on top of the other band, at which point one arrives at the original diagram of $K$.
\end{proof}

The following result was shown for the figure-eight knot (the case that the string link $J$ is a single twist) by the first author (inspired by ~\cite{FS1}). It was extended, by Cha, to the case that  $J$ is an arbitrary number of twists in ~\cite[p.63]{Cha2}. Our contribution here is just to note that Cha's proof suffices to prove this more general result.

\begin{lem}\label{lem:ratslice} Each knot $K$ in the family shown in Figure~\ref{fig:generalamphi} is slice in a rational homology $4$-ball.
\end{lem}
\begin{proof} We follow the argument of \cite{Cha2}, only indicating where our more general argument deviates. It suffices to show that the zero-framed surgery, $M_K$, as shown on the left-hand side of Figure~\ref{fig:ratball1}, is rational homology cobordant to $S^1\times S^2$. After adding, to $M_K\times [0,1]$, a four-dimensional $1$-handle and $2$-handle (going algebraically twice over the $1$-handle)  and performing certain handle slides (see \cite[p.62-64]{Cha2}), one arrives at a $3$-manifold $M'$ given by surgery on the $3$-component link drawn as the solid lines on the right-hand side of Figure~\ref{fig:ratball1}. Therefore $M_K$ is rationally homology cobordant to $M'$.

\begin{figure}[htbp]
\setlength{\unitlength}{1pt}
\begin{picture}(360,185)
\put(-20,0){\includegraphics{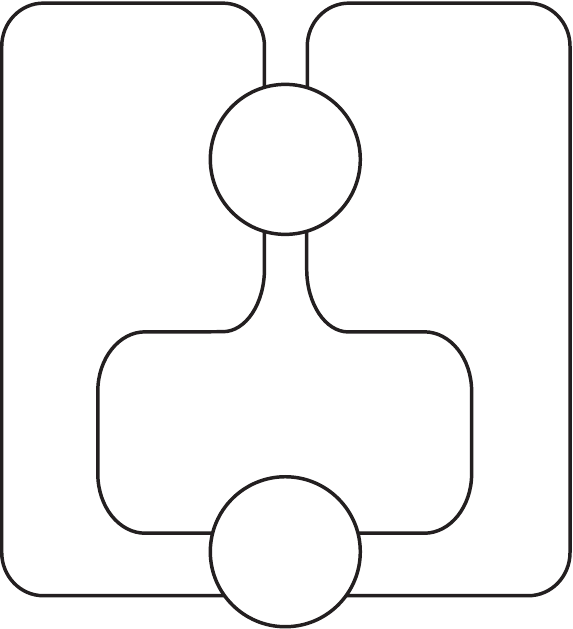}}
\put(180,0){\includegraphics{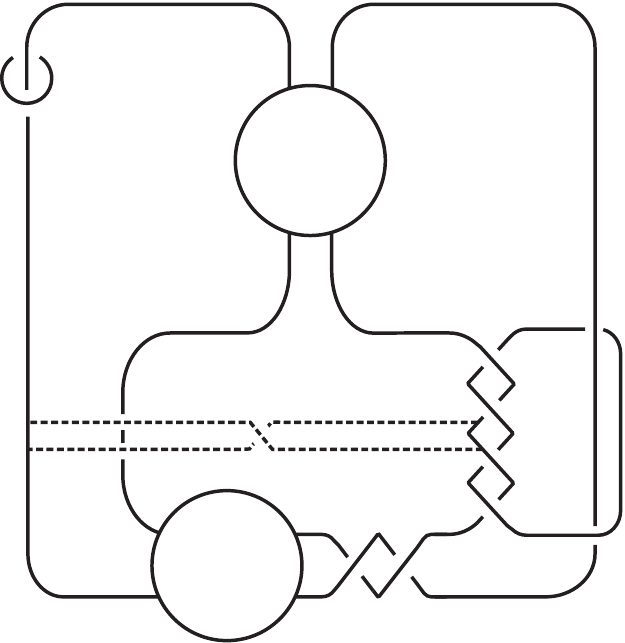}}
\put(57,134){$J$}
\put(265,136){$J$}
\put(-15,96){$0$}
\put(194,96){$0$}
\put(198,160){$2$}
\put(363,50){$-2$}
\put(57,26){\rotatebox{-90}{$\overline{J}$}}
\put(241,26){\rotatebox{-90}{$\overline{J}$}}
\end{picture}
\caption{}\label{fig:ratball1}
\end{figure}
\noindent Next one shows, as follows, that this underlying $3$-component link, $L_1$, is concordant to the simple $3$-component link, $L_4$ shown on the right-hand side of Figure~\ref{fig:ratball2}.
\begin{figure}[htbp]
\setlength{\unitlength}{1pt}
\begin{picture}(220,204)
\put(-99,0){\includegraphics{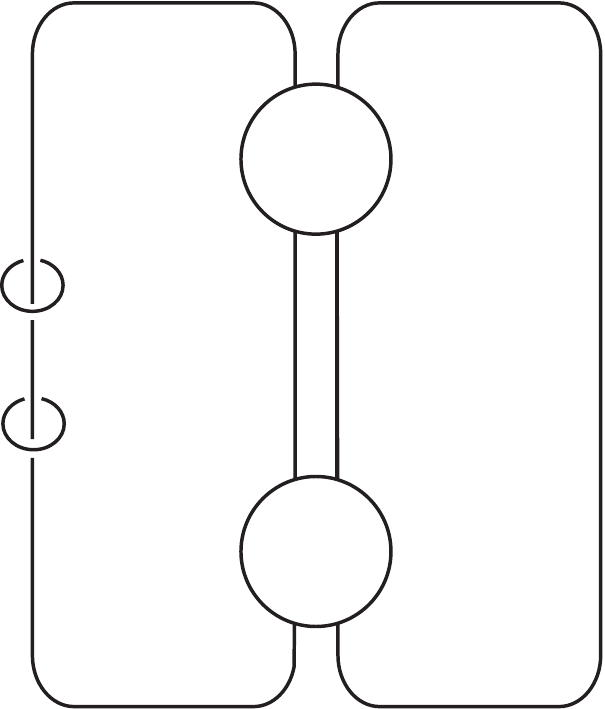}}
\put(185,50){\includegraphics{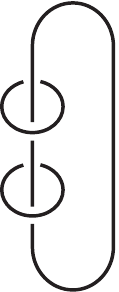}}
\put(-16,43){$-J$}
\put(-12,155){$J$}
\put(174,102){$2$}
\put(167,76){$-2$}
\put(224,65){0}
\end{picture}
\caption{}\label{fig:ratball2}
\end{figure}
Ignoring the framings on $L_1$, add a band as shown by the dashed lines on the right-hand side of Figure~\ref{fig:ratball1}, resulting in the $4$-component link, $L_2$, shown on the left-hand side of Figure~\ref{fig:ratball2}. Here $-J$ is the image of the (unoriented) string link under the map $(x,y,z)\mapsto (x,-y,z)$. One must be careful here since $\overline{J}$, which is reflection in the plane of the paper, is not the correct notion of mirror image for a string link. Since our $y$-axis is the true axis of the string link (the $[0,1]$ factor in $D^2\times [0,1]$), $-J$ is the concordance inverse of $J$ in the string link concordance group \cite{LD}\cite{HL2}, so $J+ (-J)$ is concordant to the trivial $2$-string link. Hence the link $L_2$ is concordant to the $4$-component link, $L_3$, that would result from taking $J$ to be trivial. Capping off the right-most unknotted component of $L_3$, we arrive at the $3$-component link, $L_4$, shown on the right-hand side of Figure~\ref{fig:ratball2}. This describes the desired concordance from $L_1$ to $L_4$. Consequently, $M'$ is homology cobordant to the $3$-manifold described by the framed link on  the right-hand side of Figure~\ref{fig:ratball2}, which is known to homeomorphic to $S^1\times S^2$.
\end{proof}

In this paper we will only need the special case of these lemmas wherein the string link $J$ consists of two twisted parallels of a single knotted arc as indicated by the examples in Figures~\ref{fig:Em} and ~\ref{fig:Em(J)}. Here an $m$ inside the rectangle indicates $m$ full positive twists between the two strands, and the $J$ inside the rectangle indicates that the trivial two component string link has been replaced by two parallel \emph{zero-twisted} copies of a single knotted arc $J$. This is explained more fully in Subsection~\ref{subsec:doublingoperators}.
\begin{figure}[htbp]
\setlength{\unitlength}{1pt}
\begin{picture}(171,132)
\put(0,0){\includegraphics{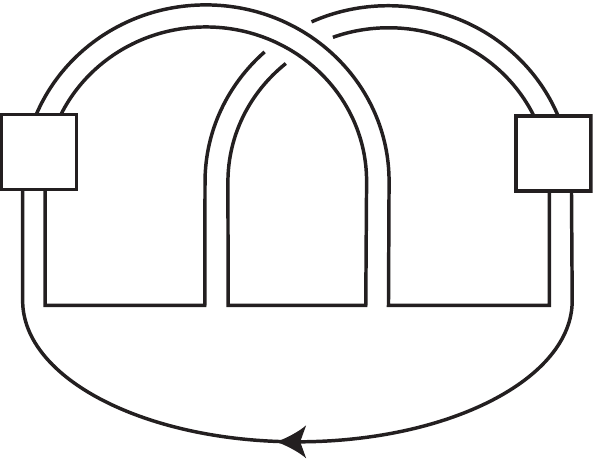}}
\put(6,85){$m$}
\put(150,84){$-m$}
\end{picture}
\caption{Negative Amphichiral Knots $\mathbb{E}^{m}$}\label{fig:Em}
\end{figure}
\begin{figure}[htbp]
\setlength{\unitlength}{1pt}
\begin{picture}(175,132)
\put(0,0){\includegraphics{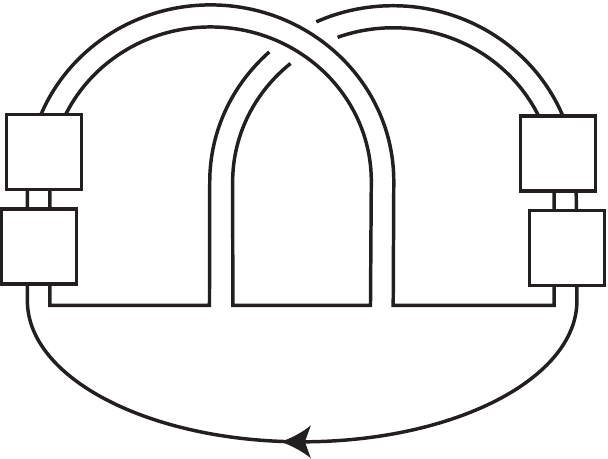}}
\put(7,85){$m$}
\put(6,57){$J$}
\put(152,85){$-m$}
\put(159,56){$\overline{J}$}
\end{picture}
\caption{Families of Negative Amphichiral Knots $\mathbb{E}^m(J)$}\label{fig:Em(J)}
\end{figure}
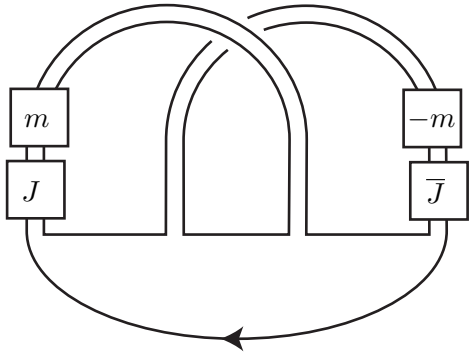

\begin{prop}\label{prop:alexpolyscoprime} If $m$ and $n$ are distinct positive integers then the Alexander polynomials $\Delta_m(t)$ of $\mathbb{E}^m$ and $\Delta_n(t)$ of $\mathbb{E}^n$ are distinct and irreducible, hence coprime. \end{prop}
\begin{proof} A Seifert matrix for $\mathbb{E}^m$ with respect to the obvious basis is
$$
\left(\begin{matrix} m & 0\cr -1 & -m
\end{matrix}\right).
$$
Thus the Alexander polynomial of $\mathbb{E}^m$ is $\Delta_m(t)=m^2t^2-(2m^2+1)t+m^2$. The discriminant $4m^2+1$ is easily seen, for $m\neq 0$, to never be the square of an integer, so the roots of $\Delta_m(t)$ are real and irrational. Hence $\Delta_m(t)$ is irreducible over $\Q[t,t^{-1}]$. It follows that if $\Delta_m(t)$ and $\Delta_n(t)$ had a common factor then they would be identical up to a unit. But the equations $m^2=qn^2$ and $2m^2+1=q(2n^2+1)$ imply $q=1$ so $m=\pm n$.
\end{proof}

\end{subsection}

\begin{subsection}{Doubling Operators}\label{subsec:doublingoperators}

To construct knots that lie deep in the $n$-solvable filtration, we use
iterated generalized satellite operations.

Suppose $R$ is a knot in $S^3$ and $\vec{\alpha}=(\alpha_1,\alpha_2,\ldots,\alpha_m)$ be an ordered, oriented, trivial link in $S^3$, that misses $R$, bounding a collection of oriented disks that meet $R$ transversely as shown on the left-hand side of Figure~\ref{fig:infection}. Suppose $(K_1,K_2,\ldots,K_m)$ is an $m$-tuple of auxiliary knots. Let $R_{\vec{\alpha}}(K_1,\ldots,K_m)$ denote the result of the operation pictured in Figure~\ref{fig:infection}, that is, for each $\alpha_j$, take the embedded disk in $S^3$ bounded by $\alpha_j$; cut off $R$ along the disk; grab the cut strands, tie them into the knot $K_j$ (such that the strands have linking number zero pairwise) and reglue as shown schematically on the right-hand side of Figure~\ref{fig:infection}.

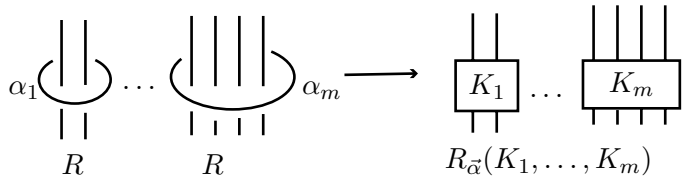
\begin{figure}[htbp]
\setlength{\unitlength}{1pt}
\begin{picture}(262,71)
\put(9,37){$\alpha_1$} \put(120,37){$\alpha_m$} \put(52,39){$\dots$}
\put(206,36){$\dots$} \put(183,37){$K_1$} \put(236,38){$K_m$}
\put(174,9){$R_{\vec{\alpha}}(K_1,\ldots,K_m)$}
\put(29,7){$R$} \put(82,7){$R$}
\put(20,20){\includegraphics{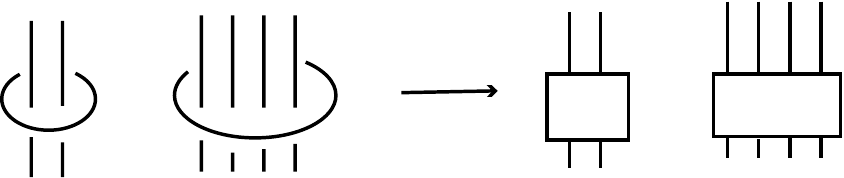}}
\end{picture}
\caption{$R_{\vec{\alpha}}(K_1,\ldots,K_m)$:
Infection of $R$ by $K_j$ along $\alpha_j$}\label{fig:infection}
\end{figure}
\noindent We will call this the \textbf{result of infection performed on the knot} $\boldsymbol{R}$ \textbf{using the infection knots} $\boldsymbol{K_j}$ \textbf{along the curves} $\boldsymbol{\alpha_j}$ ~\cite{COT2}. In the case that $m=1$ this is the same as the classical satellite construction. This construction has an alternative description. For each $\alpha_j$, remove a tubular neighborhood of $\alpha_j$ in $S^3$ and glue in the exterior of a tubular neighborhood of $K_j$ along their common boundary, which is a torus, in such a way that (the longitude of) $\alpha_j$ is identified with the meridian, $\mu_j$, of $K_j$ and the meridian of $\alpha_j$ is identified with the reverse of the longitude, $\ell_j$, of $K_j$ as suggested by Figure~\ref{fig:infectionannotated}. The resulting space can be seen to be homeomorphic to $S^3$ and the image of $R$ is the new knot.
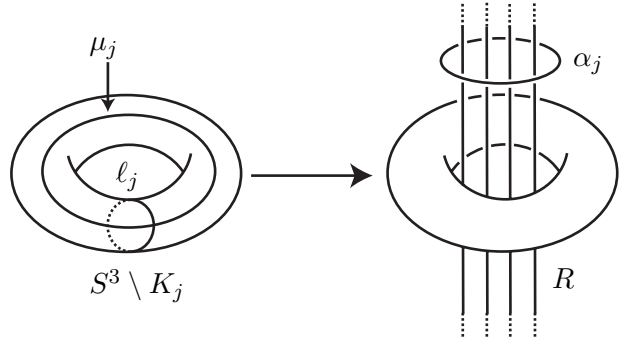
\begin{figure}[htbp]
\setlength{\unitlength}{1pt}
\begin{picture}(230,127)
\put(40,59){$\ell_j$}
\put(30,109){$\mu_j$}
\put(213,103){$\alpha_j$}
\put(205,19){$R$}
\put(30,17){$S^3\setminus K_j$}
\put(0,0){\includegraphics{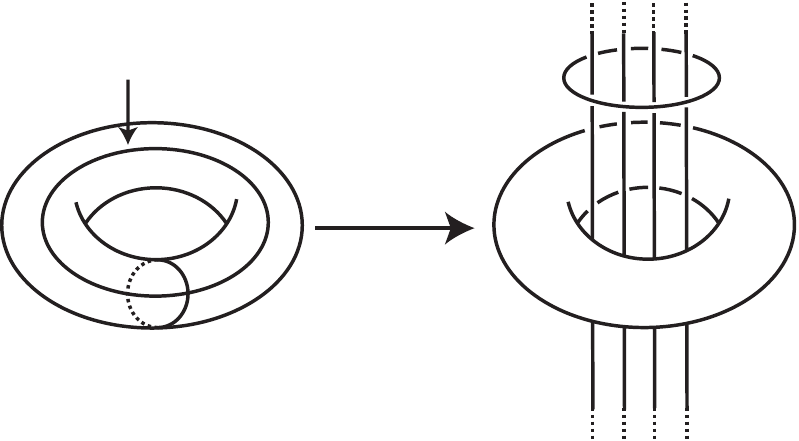}}
\end{picture}
\caption{Infection as replacing a solid torus by a knot exterior}\label{fig:infectionannotated}
\end{figure}

It is well known that if the input knots $K_{1}$ and $K_{2}$ are concordant, then the output knots $R_{\alpha}(K_{1})$ and $R_{\alpha}(K_{2})$ are concordant. Thus the functions $R_{\vec{\alpha}}$ descend to $\mathcal{C}$.

\begin{defn}[\cite{CHL4,CHL5}]\label{def:doublingop} A \textbf{doubling operator}, $R_{\vec{\alpha}}:\mathcal{C}\times\dots\times \mathcal{C}\to \mathcal{C}$ is a function, as in Figure~\ref{fig:infection}, that is given by infection on a ribbon knot $R$ wherein, for each $i$, $lk(R,\alpha_i)=0$. Often we suppress $\alpha_i$ from the notation.
\end{defn}

These are called doubling operators because they generalize untwisted Whitehead doubling.

In particular we will consider the family of doubling operators $\mathfrak{R}^{m}_{\eta_1,\eta_2}(-,-)$ shown in Figure~\ref{fig:frakR}.
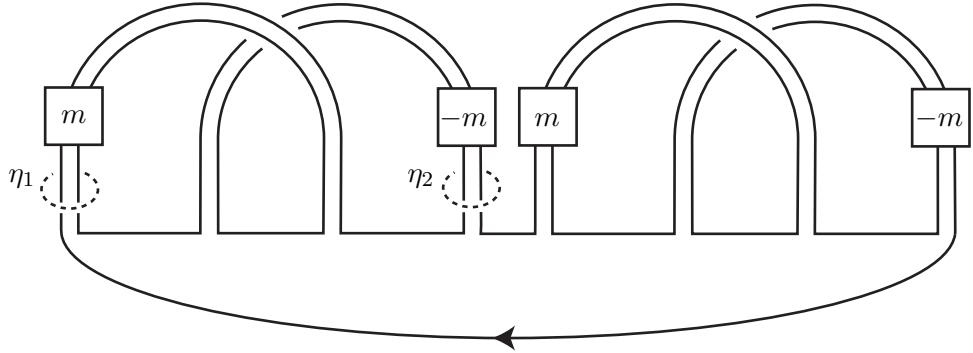
\begin{figure}[htbp]
\setlength{\unitlength}{1pt}
\begin{picture}(351,132)
\put(0,0){\includegraphics{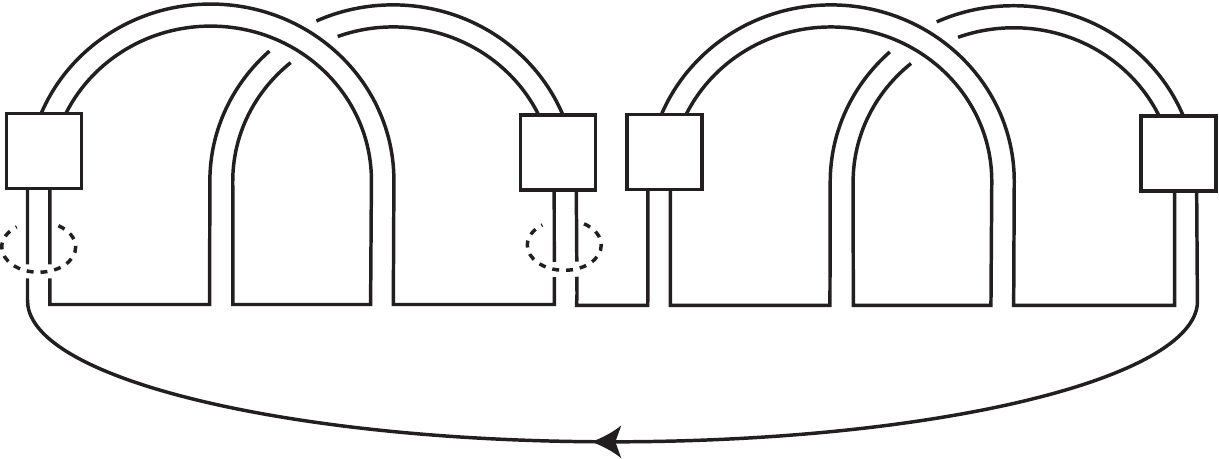}}
\put(8,86){$m$}
\put(187,85){$m$}
\put(151,85){$-m$}
\put(331,85){$-m$}
\put(139,64){$\eta_2$}
\put(-12,64){$\eta_1$}
\end{picture}
\caption{Negative Amphichiral Doubling Operators $\mathfrak{R}^m\equiv \mathbb{E}^m\# \mathbb{E}^m$}\label{fig:frakR}
\end{figure}
Note that, since $\mathbb{E}^m$ is negative amphichiral by Lemma~\ref{lem:amphi},
$$\
\mathfrak{R}^m\equiv \mathbb{E}^m\# \mathbb{E}^m\cong \mathbb{E}^m\# -\mathbb{E}^m,
$$
which is well known to be a ribbon knot \cite[Exercise 8E.30]{R}. Thus $\mathfrak{R}^{m}$ is a negative amphichiral ribbon knot. For the case $m=1$, this was already noted in ~\cite{Li6}.

We will also consider the family of doubling operators, $\R^m_\alpha$, shown in Figure~\ref{fig:ribbonfamily} (where here the $-m$ inside a box symbolizes $m$ full negative twists between the bands but where the individual bands remain untwisted), equipped with a specified circle $\alpha$ that can be shown to generate its Alexander module.

\begin{figure}[htbp]
\setlength{\unitlength}{1pt}
\begin{picture}(327,151)
\put(85,0){\includegraphics{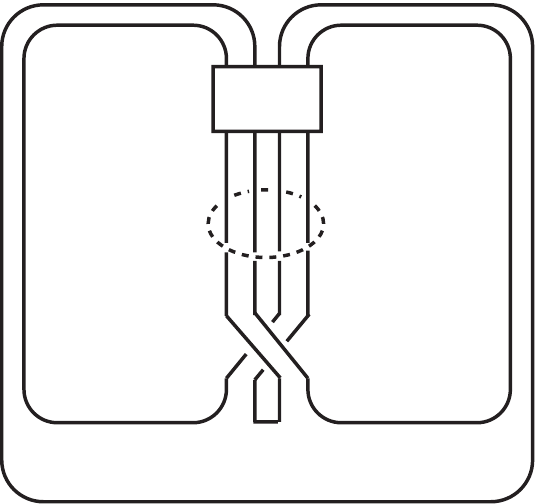}}
\put(50,78){$\R^{m}$}
\put(135,78){$\alpha$}
\put(153,114){$-m$}
\end{picture}
\caption{Doubling operators $\R^{m}_\alpha$}\label{fig:ribbonfamily}
\end{figure}

\end{subsection}

\begin{subsection}{Elements of order 2 in $\boldsymbol{\mathcal{F}_n}$}\label{subsec:ourexamples}

Now we describe large families of examples of negative amphichiral knots that lie in $\mathcal{F}_n$. Let $K^0$ be any knot with Arf($K^0$)$=0$. Let $K^{n-1}$ be the image of $K^0$ under the composition of \emph{any} $n-1$ doubling operators (each requiring a single input), that is,
\begin{equation}\label{eq:defofgeneralKn}
K^{n-1}\equiv R^{{n-1}}\circ\dots\circ R^{{1}}(K^0).
\end{equation}
Then, for any integer $m$ we define $K^n$ as in Figure~\ref{fig:examplesKn}, that is, $K^n\equiv \mathfrak{R}^{m}_{\eta_1,\eta_2}(K^{n-1},\overline{K^{n-1}})$.
\begin{figure}[htbp]
\setlength{\unitlength}{1pt}
\begin{picture}(351,132)
\put(0,0){\includegraphics{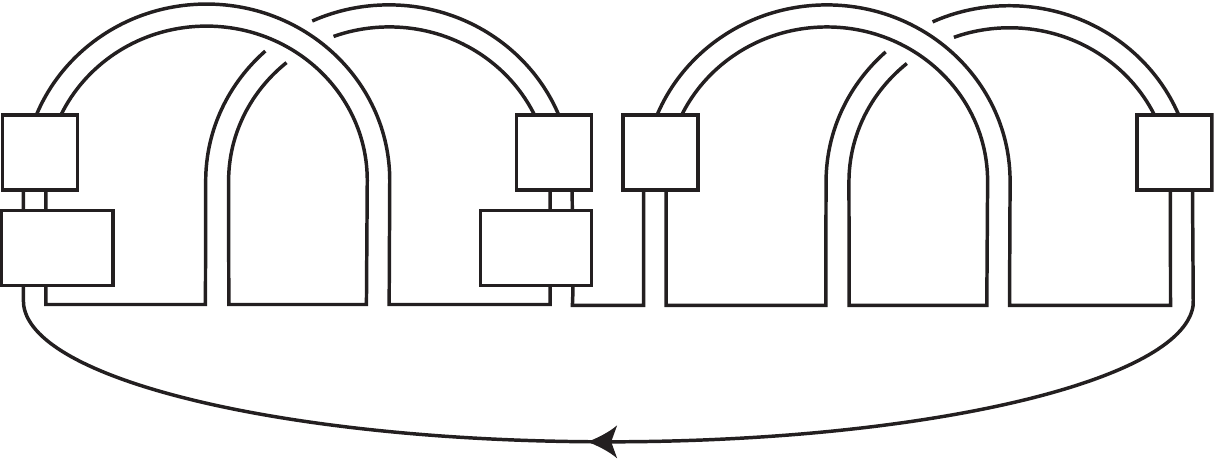}}
\put(6,86){$m$}
\put(3,57){$K^{n-1}$}
\put(186,86){$m$}
\put(142,57){$\overline{K^{n-1}}$}
\put(151,86){$-m$}
\put(330,86){$-m$}
\end{picture}
\caption{The examples $K^n$}\label{fig:examplesKn}
\end{figure}
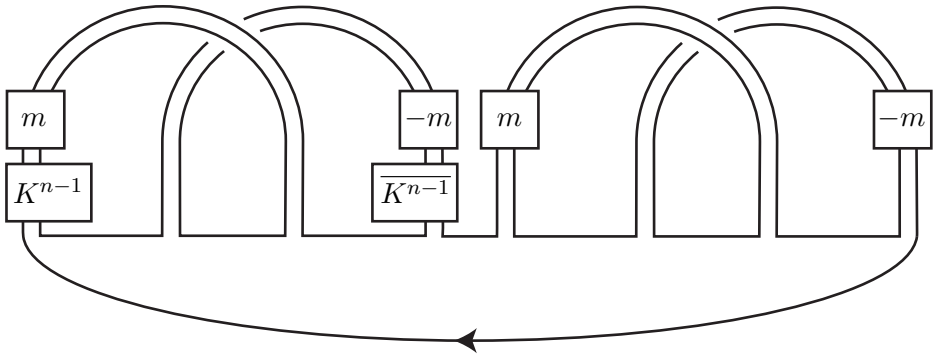

\begin{prop}\label{prop:Knisorder2} For any $n\geq 1$, any $m$, any composition of $n-1$ doubling operators and any Arf invariant zero input knot $K^0$, the knot $K^n$ of Figure~\ref{fig:examplesKn} satisfies
\begin{itemize}
\item $K^n$ is negative amphichiral;
\item $K^n\in \mathcal{F}_n$;
\item $K^n$ is (smoothly) slice in a smooth rational homology $4$-ball; and
\item $K^n\# ~K^n$ is a slice knot.
\end{itemize}
\end{prop}
\begin{proof} It was shown in \cite[Theorem 7.1]{CHL3} that, for any any doubling operator $R$,
$$
R(\mathcal{F}_i,\dots,\mathcal{F}_i)\subset \mathcal{F}_{i+1}.
$$
Since any knot of Arf invariant zero is known to lie in $\mathcal{F}_0$ \cite[Remark 8.14,Thm. 8.11]{COT}, and since $K^n$ is the image of $K^0$ under a composition of $n$ doubling operators, it follows that $K^n\in \mathcal{F}_n$.

Note that $K^n$ is the connected sum of two knots each of which is of the form shown in Figure~\ref{fig:Em(J)} (hence of the form of Figure~\ref{fig:generalamphi}).  Thus, by Lemma~\ref{lem:ratslice},  each such $K^n$ is slice in a rational homology $4$-ball. Moreover, by Lemma~\ref{lem:amphi}, $K^n$ is negative amphichiral so $K^n\# K^n$ is isotopic to $K^n\# r(\overline{K^n})$. But the latter is a ribbon knot and hence a slice knot.
\end{proof}

For specificity we define the following infinite families:

\begin{defn}\label{def:Kn(m1,mn)} Given an $n$-tuple $(m_1,...,m_{n})$ of integers and an Arf invariant zero knot $K^0$, we define $\K^n(m_1,\dots,m_n,K^0)$ to be the image of $K^0$ under the following composition of $n$ doubling operators. Specifically let
$$
\K^n\equiv \K^n(m_1,\dots,m_n,K^0)\equiv  \mathfrak{R}^{m_n}_{\eta_1,\eta_2}(\K^{n-1},\overline{\K^{n-1}}),
$$
as shown in Figure~\ref{fig:examplesKn}, where $\K^{n-1}$ is
$$
\R^{m_{n-1}}\circ\dots\circ \R^{m_{1}}(K^0),
$$
where the $\R^{m_i}$ are the operators of Figure~\ref{fig:ribbonfamily}. In other words, we recursively set:
\begin{align*}
&\K^1=\R_\alpha^{m_{1}}(K^0);\\
&\K^{2}=\R_\alpha^{m_{2}}\circ \R_\alpha^{m_{1}}(K^0);\\
&\vdots\\
&\K^{n-1}=\R_\alpha^{m_{n-1}}\circ \dots \circ \R_\alpha^{m_{1}}(K^0);\\
&\K^n= \mathfrak{R}^{m_n}_{\eta_1,\eta_2}(\K^{n-1},\overline{\K^{n-1}}).
\end{align*}
Even though $\K^n$ depends on $(m_1,...,m_{n},K^0)$, we will often suppress the latter from the notation.
\end{defn}

\end{subsection}

\section{Commutator Series and Filtrations of the knot concordance groups}\label{sec:series}

To accomplish our goals, we must establish that many of the knots in the families given by Figure~\ref{fig:examplesKn}, and specifically those in Definition~\ref{def:Kn(m1,mn)}, are not in $\mathcal{F}_{n.5}$ and, indeed, are distinct from each other in $\mathcal{F}_n/\mathcal{F}_{n.5}$. To this end we review recent work of the authors that introduced refinements of the $n$-solvable filtration parameterized by certain classes of group series that generalized the derived series. In particular the authors defined specific filtrations of $\mathcal{C}$ that depend on a sequence of polynomials. These filtrations can then be used to distinguish between knots with different Alexander modules or different higher-order Alexander modules. All of the material in this section is a review of the relevant terminology of ~\cite[Sections 2,3]{CHL5}.

Recall that the \textbf{derived series}, $\{G^{(n)}~|~n\geq 0\}$, of a group $G$ is defined recursively by $G^{(0)}\equiv G$ and $G^{(n+1)}\equiv [G^{(n)},G^{(n)}]$. The \textbf{rational derived series} \cite{Ha1}, $\{G_r^{(n)}~|~n\geq 0\}$, is defined by $G_r^{(0)}\equiv G$ and
$$
G^{(n+1)}_r=\ker \left(G_r^{(n)}\to \frac{G_r^{(n)}}{[G_r^{(n)},G_r^{(n)}]}\to \frac{G_r^{(n)}}{[G_r^{(n)},G_r^{(n)}]}\otimes_\Z \Q \right).
$$

More generally,

\begin{defn}[{\cite[Definition 2.1]{CHL5}}]\label{def:series} A \textbf{commutator series} defined on a class of groups is a function, $*$, that assigns to each group $G$ in the class a nested sequence of normal subgroups
$$
\dots\vartriangleleft\gnp _*\vartriangleleft \gn_*\vartriangleleft\dots\vartriangleleft \ensuremath{G^{\sss (0)}}_*\equiv G,
$$
such that $\gn_*/\gnp_*$ is a torsion-free abelian group.
\end{defn}

\begin{prop}[{\cite[Proposition 2.2]{CHL5}}]\label{prop:commseriesprops} For any commutator series $\{\gn_*\}$,
\begin{itemize}
\item [1.] $\{x\in\gn_*~|~\exists k>0, ~x^k\in[\gn_*,\gn_*]\}\subset \gnp_*$ (and in particular $[\gn_*,\gn_*]\subset \gnp_*$, whence the name commutator series);
\item [2.] $\gn_r\subset \gn_*$, that is, every commutator series  contains the rational derived series;
\item [3.] $G/\gn_*$ is a poly-(torsion-free abelian) group (abbreviated PTFA);
\item [4.] $\Z [G/\gn_*]$ and $\Q [G/\gn_*]$ are right (and left) Ore domains.
\end{itemize}
\end{prop}

Any commutator series that satisfies a weak functoriality condition induces a filtration, $\{\mathcal{F}^*_n\}$, of $\mathcal{C}$ by subgroups. These filtrations generalize and refine the ($n$)-solvable filtration $\{\mathcal{F}_n\}$ of ~\cite{COT}. Let $M_K$ denote the closed $3$-manifold obtained by zero-framed surgery on $S^3$ along $K$.

\begin{defn}[{\cite[Definition 2.3]{CHL5}}]\label{def:Gnsolvable} A knot $K$ is an element of $\boldsymbol{\mathcal{F}_{n}^*}$ if
the zero-framed surgery $M_K$ bounds a compact, connected, oriented, smooth $4$-manifold $W$ such that
\begin{itemize}
\item [1.] $H_1(M_K;\Z)\to H_1(W;\Z)$ is an isomorphism;
\item [2.] $H_2(W;\Z)$ has a basis consisting of connected, compact, oriented surfaces, $\{L_i,D_i~|~1\leq i\leq r\}$,  embedded in $W$ with trivial normal bundles, wherein the surfaces are pairwise disjoint except that, for each $i$, $L_i$ intersects $D_i$ transversely once with positive sign;
\item [3.] for each $i$, $\pi_1(L_i)\subset \pi_1(W)^{(n)}_*$
and $\pi_1(D_i)\subset \pi_1(W)^{(n)}_*$.
\end{itemize}
A knot $K\in \boldsymbol{\mathcal{F}_{n.5}^*}$ if in addition,
\begin{itemize}
\item [4.] for each $i$, $\pi_1(L_i)\subset \pi_1(W)^{(n+1)}_*$.
\end{itemize}
Such a $4$-manifold is called an $\boldsymbol{(n,*)}$\textbf{-solution} (respectively an $(n.5,*)$-solution) for $K$ and it is said that $K$ is $\boldsymbol{(n,*)}$-\textbf{solvable} (respectively $(n.5,*)$-solvable) via $W$. The case where the commutator series is the derived series (without the torsion-free abelian restriction) is denoted $\mathcal{F}_{n}$ and we speak of $W$ being an ($n$)-solution, and $K$ or $M_K$ being ($n$)-solvable via $W$ ~\cite[Section 8]{COT}.
\end{defn}

\begin{defn}\label{def:weakfunctorial} A commutator series $\{G^{(n)}_*\}$ is \textbf{weakly functorial} (on a class of $\{$groups, maps$\}$) if $f(G^{(n)}_*)\subset \pi^{(n)}_*$ for each $n$ and for any homomorphism $f:G\to \pi$ (in the class) that induces an isomorphism $G/G^{(1)}_r\cong \pi/\pi^{(1)}_r$ (i.e. induces an isomorphism on $H_1(-;\Q)$).
\end{defn}

\begin{prop}[{\cite[Prop. 2.5]{CHL5}}]\label{prop:functseriesfilters} Suppose $*$ is a weakly functorial commutator series defined on the class of groups with $\beta_1=1$. Then  $\{\mathcal{F}_n^*\}_{n\geq 0}$ is a filtration by subgroups of the classical (smooth) knot concordance group $\mathcal{C}$:
$$
\cdots \subset \mathcal{F}_{n+1}^* \subset \mathcal{F}_{n.5}^*\subset\mathcal{F}_n^*\subset\cdots \subset
\mathcal{F}_1^*\subset \mathcal{F}_{0.5}^* \subset \mathcal{F}_{0}^* \subset \mathcal{C}.
$$
Moreover, for any $n\in \frac{1}{2}\Z$
$$
\mathcal{F}_n\subset\mathcal{F}_n^*.
$$
\end{prop}

The case where the commutator series is the derived series (without the torsion-free abelian restriction) is the $\boldsymbol{(n)}$\textbf{-solvable filtration} ~\cite{COT}, denoted $\boldsymbol{\{\mathcal{F}_{n}\}}$.

\begin{subsection}{The Derived Series Localized at $\boldsymbol{\mathcal{P}}$}\label{subsec:derivedserieslocalatp}

\

Fix an $n$-tuple $\mathcal{P}=(p_1(t),...,p_n(t))$ of non-zero elements of $\Q[t,t^{-1}]$, such that $p_1(t)\doteq p_1(t^{-1})$. For each such $\mathcal{P}$ we now recall from \cite{CHL5} the definition of a partial commutator series that we call the \emph{(polarized) derived series localized at $\mathcal{P}$}, that is defined on the class of groups with $\beta_1=1$.

Suppose $G$ is a group such that $G/G_r^{(1)}\cong \Z=\langle \mu\rangle$. Then we define the derived series localized at $\mathcal{P}$ recursively in terms of certain right divisor sets $S_{p_n}\subset \Q[G/\gn_\mathcal{P}]$.

\begin{defn}\label{def:defderivedlocalp}  For $n\geq 0$, let
$$
\begin{array}{lcc}
G^{(0)}_\mathcal{P}\equiv G;\\
G^{(1)}_\mathcal{P}\equiv G^{(1)}_r;
\end{array}
$$
and for $n\geq 1$
\begin{eqnarray}\label{eqnarray:defderivedp}
\gnp_\mathcal{P}\equiv \ker \left(\gn_\mathcal{P}\to \frac{\gn_\mathcal{P}}{[\gn_\mathcal{P} ,\gn_\mathcal{P} ]}\otimes_{\Z [G/\gn_\mathcal{P}]}\Q[G/\gn_\mathcal{P}]S_{p_n}^{-1}\right).
\end{eqnarray}
\end{defn}

To make sense of \eqref{eqnarray:defderivedp} one must realize that, for any $H\vartriangleleft G$, $H/[H,H]$ is a right $\Z [G/H]$-module where $g$ acts on $h$ by $h\mapsto g^{-1}hg$. One must also verify, at each stage, that $\gn_\mathcal{P}$ has been defined in such a way that $G^{(k)}_\mathcal{P}/G^{(k+1)}_\mathcal{P}$ is a torsion-free abelian group for each $k<n$, so $G/\gn_\mathcal{P}$ is a poly-(torsion-free-abelian) group (PTFA), from which it follows that $\Q[G/\gn_\mathcal{P}]$ is a right Ore domain. Therefore, for any right divisor set
$S_{p_n}\subset \Q[G/\gn_\mathcal{P}]$  we may define the Ore localization  $\Q[G/\gn_\mathcal{P}]S_{p_n}^{-1}$ as in \eqref{eqnarray:defderivedp} (see \cite[Sections $3$,$4$]{CHL5}).

For the (polarized) derived series localized at $\mathcal{P}$ we use the following right divisor sets:

\begin{defn}\label{def:Sdefderivedlocalp} The \textbf{(polarized) derived series localized at} $\mathcal{P}$  is defined as in Definition~\ref{def:defderivedlocalp} by setting
\begin{eqnarray}
S_{p_1}=S_{p_1}(G)=\{ q_1(\mu)...q_r(\mu)~|~(p_1(t),q_j(t))=1; ~G/G^{(1)}_r\cong\langle\mu\rangle\}\subset \Q[G/G^{(1)}_r];
\end{eqnarray}
and for $n\geq 2$
\begin{eqnarray}
S_{p_n}=S_{p_n}(G)=\{ q_1(a_1)...q_r(a_r)~|~\widetilde{(p_n,q_j)}=1; ~q_j(1)\neq 0; ~a_j\in G^{(n-1)}_\mathcal{P}/G^{(n)}_\mathcal{P}\},
\end{eqnarray}
so $S_{p_n}\subset \Q[G^{(n-1)}_\mathcal{P}/G^{(n)}_\mathcal{P}]$.
\end{defn}

Here $p_i(t)$ and $q_j(t)$ are in $\Q[t,t^{-1}]$. By $(p_1,q_j)=1$ we mean that $p_1$ is coprime to $q_j$ in $\Q[t,t^{-1}]$, as usual. But by $~\widetilde{(p_n,q_j)}=1$ we mean something slightly stronger.

\begin{defn}[{\cite[Defn. 4.4]{CHL5}}]\label{def:stronglycoprime}  Two non-zero polynomials $p(t),q(t)\in\Q [t,t^{-1}]$ are said to be \textbf{strongly coprime}, denoted $\widetilde{(p,q)}=1$ if, for every pair of non-zero integers, $n,k$, $p(t^n)$ is relatively prime to $q(t^k)$. Alternatively, $\widetilde{(p,q)}\neq 1$ if and only if there exist \emph{non-zero} roots, $r_p, r_q\in \mathbb{C}$*, of $p(t)$ and $q(t)$ respectively, and non-zero integers $k,n$, such that $r_p^k=r_q^n$. Clearly, $\widetilde{(p,q)}= 1$ if and only if for each prime factor $p_i(t)$ of $p(t)$ and $q_j(t)$ of $q(t)$, $\widetilde{(p_i,q_j)}= 1$.
\end{defn}

Note that $\Q-\{0\}\subset S_{p_n}$ (take $q_j$ to be a non-zero constant). It is easy to see (and was proved in \cite[Section 4]{CHL5}) that $S_{p_n}$ is closed (up to units) under the involution on $\Q[G/\gn_\mathcal{P}]$. Here we need $p_1(t)\doteq p_1(t^{-1})$.

\begin{ex}\label{ex:stronglycoprime} Consider the family of quadratic polynomials
$$
\{q_m(t)=(mt-(m+1))((m+1)t-m)~|~m\in \mathbb{Z}^+\},
$$
whose roots are $\{m/(m+1),(m+1)/m\}$. The polynomial $q_m$ is the Alexander polynomial of the ribbon knot $\R^m$ shown in in Figure~\ref{fig:ribbonfamily}. It can easily be seen (and was proved in ~\cite[Example 4.10]{CHL5}) that $\widetilde{(q_m,q_n)}=1$ if $m\neq n$.
\end{ex}

\begin{thm}[{\cite[Thm. 4.16]{CHL5}}]\label{thm:derivedlocpfunctorial} The (polarized) derived series localized at $\mathcal{P}$ is a weakly functorial commutator series on the class of groups with $\beta_1=1$.
\end{thm}

\end{subsection}

\begin{section}{von Neumann signature defects as obstructions to $(n.5,*)$-solvability}\label{sec:rhoinvs}

To each commutator series there exist signature defects that offer obstructions to a given knot lying in a term of $\mathcal{F}^{*}$. Given a closed, oriented 3-manifold $M$, a discrete group $\G$, and a representation $\phi : \pi_1(M)
\to \G$, the \textbf{von Neumann
$\boldsymbol{\rho}$-invariant}, $\rho(M,\phi)\in \mathbb{R}$, was defined by Cheeger and Gromov  ~\cite{ChGr1}. If $(M,\phi) = \partial
(W,\psi)$ for some compact, oriented 4-manifold $W$ and $\psi : \pi_1(W) \to \G$, then it is known that $\rho(M,\phi) =
\sigma^{(2)}_\G(W,\psi) - \sigma(W)$ where $\sigma^{(2)}_\G(W,\psi)$ is the
$L^{(2)}$-signature (von Neumann signature) of the equivariant intersection form defined on
$H_2(W;\mathbb{Z}\G)$ twisted by $\psi$, and $\sigma(W)$ is the ordinary
signature of $W$ ~\cite{LS}\cite[Section 2]{CT}. Thus the $\rho$-invariants should be thought of as \emph{signature defects}. They were first used to detect non-slice knots in ~\cite{COT}. For a more thorough discussion see ~\cite[Section 5]{COT}\cite[Section 2]{CT}\cite[Section 2]{COT2}. All of the coefficient systems $\G$ in this paper will be of the form $\pi/\pi^{(n)}_*$ where $\pi$ is the fundamental group of a space. Hence all such $\G$ will be PTFA. Aside from the definition, the properties that we use in this paper are:
\begin{prop}\label{prop:rhoprops}\

\begin{itemize}
\item [1.] If $\phi$ factors through $\phi': \pi_1(M)\to \G'$ where
$\G'$ is a subgroup of $\G$, then $\rho(M,\phi') = \rho(M,\phi)$.
\item [2.] If $\phi$ is trivial (the zero map), then $\rho(M,\phi) = 0$.
\item [3.] If $M=M_K$ is the zero-surgery on a knot $K$ and $\phi:\pi_1(M)\to \mathbb{Z}$ is the abelianization, then $\rho(M,\phi)$ is denoted $\boldsymbol{\rho_0(K)}$ and is equal to the integral over the circle of the Levine-Tristram signature function of $K$ ~\cite[Prop. 5.1]{COT2}. Thus $\rho_0(K)$ is the average of the classical signatures of $K$.
\item [4.] If $K$ is a slice knot or link and $\phi:M_K\to \G$ ($\G$ PTFA) extends over $\pi_1$ of a slice disk exterior then $\rho(M_K,\phi)=0$ by ~\cite[Theorem 4.2]{COT}.
\item [5.] The von Neumann signature satisfies Novikov  additivity, i.e. if $W_1$ and $W_2$ intersect along a common boundary component then $\sigma^{(2)}_\G(W_1\cup W_2)=\sigma^{(2)}_\G(W_1)+\sigma^{(2)}_\G(W_2)$ ~\cite[Lemma 5.9]{COT}.
\item  [6.] For any $3$-manifold $M$, there is a positive real number $C_M$, called the \textbf{Cheeger-Gromov constant} \cite{ChGr1}\cite[Section 2]{CT} of $M$ such that, for \emph{any} $\phi$
    $$
    |\rho(M,\phi)|<C_M.
    $$
    \end{itemize}
\end{prop}

We will also need the following generalization of property ($4$).

\begin{thm}[{\cite[Theorem 5.2]{CHL5}}]\label{thm:generalsignaturesobstruct} Suppose $*$ is a commutator series (no functoriality is required). Suppose $K\in \mathcal{F}_{n.5}^*$, so
the zero-framed surgery $M_K$ is $(n.5,*)$-solvable via $W$ as in Definition~\ref{def:Gnsolvable}. Let $G=\pi_1(W)$ and consider
$$
\phi:\pi_1(M_K)\to G\to G/\gnp_*\to \G,
$$
where $\G$ is an arbitrary PTFA group.
Then
$$
\sigma^{(2)}_{\Gamma}(W,\phi)-\sigma(W)=0=\rho(M_K,\phi).
$$
\end{thm}

\end{section}

\section{Statements of Main Results and the outline of the proof}\label{sec:outline}

We will show that for any $n\geq 2$, not only does there exist
$$
\Z_2^\infty\subset \mathbb{G}_{n}\equiv \mathcal{F}_{n}/\mathcal{F}_{n.5},
$$
but there are also many distinct such classes
$$
\bigoplus_{\substack{\mathbb{P}_{n-1}}}\Z_2^\infty\subset \mathbb{G}_n,
$$
distinguished by the sequence of orders of certain higher-order Alexander modules of the knots.

Given the sequence $\mathcal{P}=(p_1(t),...,p_n(t))$, we have defined (in Definitions~\ref{def:defderivedlocalp} and ~\ref{def:Sdefderivedlocalp}) an associated commutator series called the derived series localized at $\mathcal{P}$.

\begin{defn}\label{def:Fp} Let $\{\mathcal{F}_n^\mathcal{P}\}$ denote the filtration of $\mathcal{C}$ associated, by Definition~\ref{def:Gnsolvable}, to the derived series localized at $\mathcal{P}$.
\end{defn}

\noindent Since for any group $G$ and integer $n$ (or half-integer), $G^{(n)}\subset G^{(n)}_\mathcal{P}$, one sees that $\mathcal{F}_n\subset\mathcal{F}_n^\mathcal{P}$. In particular $\mathcal{F}_{n.5}\subset\mathcal{F}_{n.5}^\mathcal{P}$, so there is a surjection
$$
\frac{\mathcal{F}_n}{\mathcal{F}_{n.5}}\overset{\pi}\twoheadrightarrow \frac{\mathcal{F}_n}{\mathcal{F}_{n.5}^\mathcal{P}}.
$$
The point of the filtration $\{\mathcal{F}_n^\mathcal{P}\}$, is that any knot $K\in \mathcal{F}_n$, whose classical Alexander polynomial is coprime to $p_1(t)$, will lie in the \textbf{kernel} of $\pi$. Moreover, the idea of Theorem~\ref{thm:KnliesinQplusone} below is that a knot will of necessity lie in the kernel of $\pi$, \emph{unless} $p_1(t)$ divides its classical Alexander polynomial \emph{and}, loosely speaking, the higher $p_i(t)$ are related to torsion in its $i^{th}$ higher-order Alexander module.

\begin{defn}\label{def:distinctP} Given $\mathcal{P}=(p_1(t),...,p_n(t))$ and $\mathcal{Q}=(q_1(t),...,q_n(t))$, we say that $\mathcal{P}$ is \textbf{strongly coprime} to $\mathcal{Q}$ if either $(p_1,q_1)=1$, or for some $k>1$, $\widetilde{(p_k,q_{k})}=1$.
\end{defn}

\begin{thm}[{\cite[Theorem 6.5]{CHL5}}]\label{thm:KnliesinQplusone} For any $n\geq 1$, let $R_{\alpha_{n-1}}^{n-1},\dots,R_{\alpha_1}^1$ be any doubling operators and $K^0$ be any Arf invariant zero input knot. Consider the knot $K^n\equiv \mathfrak{R}^{m}_{\eta_1,\eta_2}(K^{n-1},\overline{K^{n-1}})$,
where $K^{n-1}=R_{\alpha_{n-1}}^{n-1}\circ\dots\circ R_{\alpha_1}^1(K^0)$. Then
$$
K^n\in \mathcal{F}_{n+1}^\mathcal{P}
$$
for each $\mathcal{P}=(p_1(t),p_2(t),...,p_n(t))$, with $p_1(t)\doteq p_1(t^{-1})$, that is strongly coprime to \linebreak
$(\Delta_m(t),q_{n-1}(t),\dots,q_1(t))$, where $\Delta_m$ is the Alexander polynomial of $\mathbb{E}^m$ and $q_i$ is the Alexander polynomial of $R^i$.
\end{thm}

This applies, in particular, to the families of Definition~\ref{def:Kn(m1,mn)}, constructed using the ribbon knots of Figures~\ref{fig:frakR} and \ref{fig:ribbonfamily}.

\begin{cor}\label{cor:KnliesinQplusone} For any $(m_1,\dots,m_n)$ and any input knot $K^0$ with Arf invariant zero,
$$
\K^n(m_1,\dots,m_n,K^0)\in \mathcal{F}_{n+1}^\mathcal{P}
$$
for each $\mathcal{P}=(p_1(t),p_2(t),...,p_n(t))$ that is strongly coprime to $(\Delta_{m_n}(t),q_{n-1}(t)\dots,q_1(t))$ where $\Delta_{m_n}$ is the Alexander polynomial of $\mathbb{E}^{m_n}$ and $q_i$ is the Alexander polynomial of $\R^{m_i}$.
\end{cor}

Now we need a non-triviality theorem to complement Theorem~\ref{thm:KnliesinQplusone}.

\begin{thm}\label{thm:main} Suppose
$$
K^n\equiv \mathfrak{R}^{m}_{\eta_1,\eta_2}(K^{n-1},\overline{K^{n-1}}),
$$
where $K^{n-1}$ is the result of applying any sequence of $n-1$ doubling operators, $R_{\alpha_{n-1}}^{n-1}\circ\dots\circ R_{\alpha_1}^1$ to an Arf invariant zero ``input'' knot $K^0$. Suppose additionally that $n\geq 2$ and
\begin{itemize}
\item [1.] $m\neq 0$;
\item [2.] for each $i$, $\alpha_i$ generates the rational Alexander module of $R^i$, and this module is non-trivial;
\item [3.] $|\rho_0(K^0)|$, the average Levine-Tristram signature of $K^0$, is greater than twice the sum of the Cheeger-Gromov constants of the ribbon knots $\mathfrak{R}^{m}$, $R^{1},\dots, R^{n-1}$ (see Section~\ref{sec:rhoinvs}).
\end{itemize}
If $\mathcal{P}$ is the sequence of classical Alexander polynomials of the knots $(\mathbb{E}^m, R^{n-1},\dots, R^1)$, then
$$
K^n\notin \mathcal{F}_{n.5}^{\sP}.
$$
\end{thm}

This can  be applied to the specific families of Definition~\ref{def:Kn(m1,mn)}.
\begin{cor}\label{cor:main}
Fix $n\geq 2$ and an $n$-tuple of positive integers $(m_1,\dots,m_n)$. Suppose $K^0$ is chosen so that $|\rho_0(K^0)|$ is greater than twice the sum of the Cheeger-Gromov constants of the ribbon knots $\mathfrak{R}^{m_n}$, $\R^{m_{n-1}},\dots, \R^{m_1}$. If $\mathcal{P}$ is the $n$-tuple of Alexander polynomials of the knots $(\mathbb{E}^{m_n},\R^{m_{n-1}},\dots,\R^{m_{1}})$, then
$$
\K^n\notin \mathcal{F}_{n.5}^{\sP}.
$$
\end{cor}

The proofs of Theorems~\ref{thm:KnliesinQplusone} and ~\ref{thm:main}  will constitute Sections~\ref{sec:steps23} and \ref{sec:proof}.
Assuming these theorems, we now derive our main results.

\begin{thm}\label{thm:independence} Fix $n\geq 2$. For any $n$-tuple of positive integers $(m_1,...,m_n)$ choose an Arf invariant zero knot $K^0(m_1,...,m_n)$ such that $|\rho_0(K^0)|$ is greater than twice the sum of the Cheeger-Gromov constants of $\mathfrak{R}^{m_n}$, $\R^{m_{n-1}},\dots, \R^{m_1}$. Then the resulting set of knots
$$
\{\K^n(m_1,...,m_n,K^0)~~|~ m_i\in \mathbb{Z}^+\},
$$
as in Definition~\ref{def:Kn(m1,mn)}, represent linearly independent, order two elements of $\mathcal{F}_n/\mathcal{F}_{n.5}$. They also represent linearly independent order two elements in $\mathcal{C}$. In particular this gives
$$
\Z_2^\infty\subset \mathbb{G}_n \equiv \frac{\mathcal{F}_n}{\mathcal{F}_{n.5}},
$$
where each class is represented by a negative amphichiral knot that is slice in a rational homology $4$-ball.
\end{thm}

\begin{proof}[Proof of Theorem~\ref{thm:independence} assuming Theorems~\ref{thm:KnliesinQplusone} and ~\ref{thm:main}] By Proposition~\ref{prop:Knisorder2}, $\K^n$ is negative amphichiral, $\K^n\in\mathcal{F}_n$ and $\K^n\# \K^n$ is a slice knot. Thus $2[\K^n]=0$ in $\mathcal{F}_n/\mathcal{F}_{n.5}$. By Corollary~\ref{cor:main}, for a certain $\mathcal{P}$, $\K^n\notin \mathcal{F}_{n.5}^\mathcal{P}$, so in particular $\K^n\notin \mathcal{F}_{n.5}$ by Proposition~\ref{prop:functseriesfilters}. Therefore each $[\K^n]$ has order precisely two in $\mathbb{G}_n$.

Suppose there exists a nontrivial relation
$$
J=\K^n(m_{11},...,m_{1n},K^0_1)\#...\# \K^n(m_{k1},...,m_{kn},K^0_k)\in \mathcal{F}_{n.5}.
$$
Set $\mathcal{P}=(p_1,...,p_n)=(\Delta_{1n},q_{1n-1},...,q_{11})$, the reverse of the sequence of Alexander polynomials of the operators corresponding to the \emph{first summand} of $J$.  For each of the \emph{other} summands of $J$, the corresponding $n$-tuple $(m_{i1},...,m_{in})$ is assumed distinct from $(m_{11},...,m_{1n})$. Therefore, the (reversed) sequence of Alexander polynomials of the operators corresponding to this other summand is strongly coprime to $\mathcal{P}$ by Proposition~\ref{prop:alexpolyscoprime} and Example~\ref{ex:stronglycoprime}. Thus, by Theorem~\ref{thm:KnliesinQplusone}, each summand of $J$, aside from the first, lies in $\mathcal{F}_{n+1}^\mathcal{P}$ and hence in $\mathcal{F}_{n.5}^\mathcal{P}$. Since $J\in \mathcal{F}_{n.5}$, $J\in \mathcal{F}_{n.5}^\mathcal{P}$, by Proposition~\ref{prop:functseriesfilters}. Since $\mathcal{F}_{n.5}^\mathcal{P}$ is a subgroup, it would follow that the first summand of $J$ also lay in $\mathcal{F}_{n.5}^\mathcal{P}$, contradicting Corollary~\ref{cor:main}.
\end{proof}

More generally,
\begin{thm}\label{thm:bigsum} Suppose $n\geq 2$. Let $\mathbb{P}_n$ be any set of $n$-tuples $\mathcal{P}=(\delta_1(t),\delta_2(t),\dots,\delta_n(t))$ of prime polynomials $\delta_i(t)\in \Z[t,t^{-1}]$ such that $\delta_i(1)=\pm 1$, $\delta_1(t)=\Delta_m=m^2t^2-(2m^2+1)t+m^2$ and with the property that, for any distinct $\mathcal{P}$, $\mathcal{P}'\in \mathbb{P}_n$, $\mathcal{P}$
and $\mathcal{P}'$ are strongly coprime. Then there exists a set of negative amphichiral $n$-solvable knots indexed by $\mathbb{P}_n$ that is linearly independent modulo $\mathcal{F}_{n.5}$, that is, that spans
$$
\bigoplus_{\substack{\mathbb{P}_{n}}}\Z_2\subset \mathbb{G}_n,
$$
where the knot corresponding to the sequence $(\delta_1(t),\delta_2(t),\dots,\delta_n(t))$ admits a sequence of higher-order Alexander modules containing submodules whose orders are determined by the sequence $(\delta_1(t)\delta_1(t^{-1}),\dots,\delta_n(t)\delta_n(t^{-1}))$ with the classical Alexander polynomial being $\delta_1(t)\delta_1(t^{-1})$. Moreover each class is represented by a negative amphichiral knot that is slice in a rational homology $4$-ball.
\end{thm}
\begin{proof}[Proof of Theorem~\ref{thm:bigsum} assuming Theorems~\ref{thm:KnliesinQplusone} and ~\ref{thm:main}] By \cite{Ter}, for any prime $\delta(t)$ with $\delta(1)=\pm 1$ there exists a ribbon knot whose Alexander module is cyclic of order $\delta(t)\delta(t^{-1})$. Hence, given $\mathcal{P}=(\delta_1(t),\delta_2(t),\dots,\delta_n(t))$, choose such ribbon knots $R^{n-1},\dots,R^1$ whose Alexander polynomials are
$\delta_2(t)\delta_2(t^{-1}),\dots,\delta_n(t)\delta_n(t^{-1})$ respectively, and choose curves $\alpha_i$ (unknotted in $S^3$), that generate the Alexander modules of the $R^i$. Thus doubling operators $R^i_{\alpha_i}$, $1\leq i\leq n-1$, are defined. Since $\delta_1(t)=\Delta_m=m^2t^2-(2m^2+1)t+m^2$, there is a ribbon knot, namely $\mathfrak{R}^m=\mathbb{E}^m\#\mathbb{E}^m$ of Figure~\ref{fig:frakR}, whose Alexander polynomial is $\delta_1(t)\delta_1(t^{-1})$. The hypotheses imply $m\neq 0$. Choose any Arf invariant zero knot $K^0$ such that $|\rho_0(K^0)|$ is greater than twice the sum of the Cheeger-Gromov constants of $\mathfrak{R}^{m}$, $R^{n-1},\dots, R^{1}$. Then set
\begin{equation}\label{eq:mostgeneralexamples}
K^n_\mathcal{P}\equiv\mathfrak{R}^{m}_{\eta_1,\eta_2}(K^{n-1},\overline{K^{n-1}}),
\end{equation}
where $K^{n-1}\equiv R_{\alpha_{n-1}}^{n-1}\circ\dots\circ R_{\alpha_1}^1(K^0)$. To each $\mathcal{P}$ there is an associated $n$-tuple, $\mathcal{P}^*=(\delta_1,\delta_2(t)\delta_2(t^{-1}),\dots,\delta_n(t)\delta_n(t^{-1}))$, that gives the sequence of Alexander polynomials of the knots $\mathbb{E}^m,R^{n-1},\dots,R^1$ that define $K^n_\mathcal{P}$.

By Lemma~\ref{lem:amphi} and Proposition~\ref{prop:Knisorder2}, each $K^n_\mathcal{P}$ is negative amphichiral and $n$-solvable.  By Theorem~\ref{thm:main},
\begin{equation}\label{eq:Knnotin}
K^n_\mathcal{P}\notin \mathcal{F}_{n.5}^{\mathcal{P}^*},
\end{equation}
so $K^n_\mathcal{P}\notin \mathcal{F}_{n.5}$. Thus $[K^n_\mathcal{P}]$ has order precisely two in $\mathbb{G}_n$. Suppose there were a non-trivial relation
$$
J=\sum_{\substack{i=1}}^{\substack{k}}K^n_{\mathcal{P}_i}\in \mathcal{F}_{n.5}.
$$
By hypothesis, if $i\neq 1$ then $\mathcal{P}_i$ is strongly coprime to $\mathcal{P}_1$. It follows that $\mathcal{P}_i^*$ is strongly coprime to $\mathcal{P}_1^*$. Thus, by Theorem~\ref{thm:KnliesinQplusone}, if $i\neq 1$ then
$$
K^n_{\mathcal{P}_i}\in \mathcal{F}_{n+1}^{\mathcal{P}_1^*}\subset \mathcal{F}_{n.5}^{\mathcal{P}_1^*}.
$$
Since $J\in \mathcal{F}_{n.5}$, $J\in \mathcal{F}_{n.5}^{\mathcal{P}_1^*}$. Since the latter is a subgroup,
$$
K^n_{\mathcal{P}_1}\in \mathcal{F}_{n.5}^{\mathcal{P}_1^*},
$$
contradicting \eqref{eq:Knnotin}.

It remains only to relate the sequence $\mathcal{P}$ to the higher-order Alexander modules of the knots $K^n_\mathcal{P}$. Since this is not central to our results, we sketch the proof. Recall:
\begin{defn}[{\cite[Def. 2.8]{C}\cite[Def. 5.3]{Ha1}}]\label{def:highordermodule} The $i^{th}$, $i\geq 1$, \textbf{higher-order} (integral) \textbf{Alexander module of a knot} $K$ is
$$
\mathcal{A}^\Z_i(K)\equiv H_1(M_{K};\Z[G/G^{(i+1)}_r])\cong \frac{G^{(i+1)}_r}{[G^{(i+1)}_r,G^{(i+1)}_r]},
$$
where $G\equiv \pi_1(M_{K})$. Note: The case $i=0$ would give the classical Alexander module.
\end{defn}

Thus $\mathcal{A}^\Z_i(K^n_\mathcal{P})$ is a module over $\G_i\equiv G/G^{(i+1)}_r$, where $G\equiv \pi_1(M_{K^n_\mathcal{P}})$. The following lemma shows that the (two) images of the \emph{classical} Alexander polynomial, $\delta_{i+1}(t)\delta_{i+1}(t^{-1})$, of the constituent operator $R^{n-i}$ under certain maps
$$
\Z[t,t^{-1}]\to \Z[G^{(i)}/G^{(i+1)}_r]\subset \Z\G_i,
$$
wherein $t\mapsto x_1$ and  $t\mapsto x_2$, appear as the orders of cyclic submodules of $\mathcal{A}^\Z_i(K^n_\mathcal{P})$.

\begin{lem}\label{lem:Pequalsgihigheralexpoly} Fixing $\mathcal{P}=(\delta_1(t),\delta_2(t),\dots,\delta_n(t))$, for each $1\leq i\leq n-1$, the $i^{th}$ higher-order Alexander module of $K^n_\mathcal{P}$ (the knot defined in \eqref{eq:mostgeneralexamples}) contains two non-trivial summands
$$
\frac{\Z\G_i}{\delta_{i+1}(x_1)\delta_{i+1}(x_1^{-1})\Z\G_i}\oplus \frac{\Z\G_i}{\delta_{i+1}(x_2)\delta_{i+1}(x_2^{-1})\Z\G_i}
$$
for certain $x_1,x_2\in G^{(i)}/G^{(i+1)}_r$.
\end{lem}
\begin{proof} Recall that $K^n_\mathcal{P}$ is defined as the image of $K^0$ under a composition of $n$ doubling operators. In particular $K^{n-1}\equiv R_{\alpha_{n-1}}^{n-1}\circ\dots\circ R_{\alpha_1}^1(K^0)$. Sequences of satellite operations have a certain associativity property yielding, for each $i\geq 2$, an alternative description of $K^{n-1}$ as a \emph{single} infection on single ribbon knot, $\tilde{R}^i$, along a curve lying in $\pi_1(S^3-\tilde{R}^i)^{(i-1)}$, using the knot $K^{n-i}$ \cite[Prop. 4.7]{CHL1A}\cite[Prop. 5.10]{CHL4}. Specifically,
$$
\begin{array}{llc}
K^{n-1}=R_{\alpha_{n-1}}^{n-1}\circ\dots\circ R_{\alpha_{n-i+1}}^{n-i+1}\left(R^{n-i}_{\alpha_{n-i}}\dots\circ R_{\alpha_1}^1(K^0)\right)\\
K^{n-1}=R_{\alpha_{n-1}}^{n-1}\circ\dots\circ R_{\alpha_{n-i+1}}^{n-i+1}(K^{n-i})\\
K^{n-1}=\left(R_{\alpha_{n-1}}^{n-1}\circ\dots\circ R_{\alpha_{n-i+2}}^{n-i+2}(R_{\alpha_{n-i+1}}^{n-i+1})\right)_{\beta_i}\left(K^{n-i}\right)\\
K^{n-1}=\tilde{R}^i_{\beta_{i}}(K^{n-i}),
\end{array}
$$
where
$$
\tilde{R}^i_{\beta_{i}}\equiv R_{\alpha_{n-1}}^{n-1}\circ\dots\circ R_{\alpha_{n-i+2}}^{n-i+2}(R_{\alpha_{n-i+1}}^{n-i+1})
$$
and $\beta_i$ is the image of $\alpha_{n-i+1}$. The specific nature of $\tilde{R}^i$ is not important to our present considerations. If $i=1$, let $\tilde{R}^i_{\beta_{i}}$ be the identity operator. Then, for any $i\geq 1$, it follows that
$$
K^n=\mathfrak{R}^{m}_{\eta_1,\eta_2}\left(\tilde{R}^i_{\beta_{i}}(K^{n-i}),\overline{\tilde{R}^i_{\beta_{i}}(K^{n-i})}\right).
$$
This can be reformulated, by the same considerations as above, to yield
$$
K^n=\tilde{\mathfrak{R}}_{\gamma_1,\gamma_2}\left(K^{n-i},\overline{K^{n-i}}\right)
$$
where $\tilde{\mathfrak{R}}=\mathfrak{R}^{m}_{\eta_1,\eta_2}(\tilde{R}^i,\tilde{R}^i)$
and $\{\gamma_1,\gamma_2\}$ are the images of the two copies of $\beta_i$. These curves can inductively shown to lie in $\pi_1(S^3-\tilde{\mathfrak{R}})^{(i)}$ \cite[Prop. 4.7]{CHL1A}\cite[Prop. 5.10]{CHL4}. The latter computation is very similar to the computation we will perform in \eqref{eq:munminusone1}.

Now we can apply known results about the effect of single infection on the higher-order Alexander modules \cite[Theorem 3.5]{Lei3}\cite[Theorem 8.2]{C}:
$$
\mathcal{A}^\Z_i(K^n)=\mathcal{A}^\Z_i(\tilde{\mathfrak{R}})\oplus \left( \mathcal{A}_0^\Z(K^{n-i})\otimes_{\Z[t,t^{-1}]} \Z\G_i\right) \oplus \left(\mathcal{A}_0^\Z(\overline{K}^{n-i})\otimes_{\Z[t,t^{-1}]} \Z\G_i\right).
$$
where $\mathcal{A}_0^\Z$ denotes the classical Alexander module and the first tensor product is given by $t\mapsto x_1=\gamma_1$ and the second by $t\mapsto x_2=\gamma_2$. But
$$
\mathcal{A}_0^\Z(K^{n-i})\cong \mathcal{A}_0^\Z(R^{n-i})\cong \frac{\Z[t,t^{-1}]}{\delta_{i+1}(t)\delta_{i+1}^{-1}(t)\Z[t,t^{-1}]}.
$$
where $t\mapsto x_1$. The Alexander modules of $\overline{R^{n-i}}$ and $R^{n-i}$ are isomorphic.  Thus $\mathcal{A}^\Z_i(K^n)$ contains two cyclic summands as claimed. By \cite[Theorem 3.5]{Lei3}\cite[Theorem 8.2]{C} these summands are non-zero precisely when $x_1$ and $x_2$ are not zero in $\G_i$. The verification of the latter requires further computation as in \cite[Theorem 4.11]{CHL1A}\cite[Propoposition 5.14]{CHL4}. These calculations are entirely similar to and easier than the ones we will do to verify our Proposition~\ref{prop:inductivestep}. They are not included.

This concludes what we will say about the connections between $\mathcal{P}$ and the orders of the higher-order Alexander modules of $K^n_\mathcal{P}$.
\end{proof}

This concludes the proof of Theorem~\ref{thm:bigsum}.

\end{proof}

\section{Sketch of Proof of Theorem~\ref{thm:KnliesinQplusone}}\label{sec:steps23}

Theorem~\ref{thm:KnliesinQplusone} is a consequence of \cite[Theorem 6.5]{CHL5}. However, we shall sketch the proof since the basic idea is elementary and it also shows that $K^n\in \mathcal{F}_n$.

\newtheorem*{thm:KnliesinQplusone}{Theorem~\ref{thm:KnliesinQplusone}}
\begin{thm:KnliesinQplusone}[{\cite[Theorem 6.5]{CHL5}}] For any $n\geq 1$ and $m\in \Z$, let $R_{\alpha_{n-1}}^{n-1},\dots,R_{\alpha_1}^1$ be any doubling operators and $K^0$ be any Arf invariant zero input knot. Consider the knot $K^n\equiv \mathfrak{R}^{m}_{\eta_1,\eta_2}(K^{n-1},\overline{K^{n-1}})$,
where $K^{n-1}=R_{\alpha_{n-1}}^{n-1}\circ\dots\circ R_{\alpha_1}^1(K^0)$. Then
$$
K^n\in \mathcal{F}_{n+1}^\mathcal{P}
$$
for each $\mathcal{P}=(p_1(t),p_2(t),...,p_n(t))$, with $p_1(t)\doteq p_1(t^{-1})$, that is strongly coprime to \linebreak
$(\Delta_m(t),q_{n-1}(t),\dots,q_1(t))$, where $\Delta_m$ is the Alexander polynomial of $\mathbb{E}^m$ and $q_i$ is the Alexander polynomial of $R^i$.
\end{thm:KnliesinQplusone}

\begin{proof}[Proof of Theorem~\ref{thm:KnliesinQplusone}] We set $K^1=R^1(K^0),\dots,K^{i}=R^i(K^{i-1})$ for $i=1,\dots,n-1$ and $K^n=\mathfrak{R}^{m}(K^{n-1},\overline{K^{n-1}})$. Recall from ~\cite[Lemma 2.3, Figure 2.1]{CHL3} that, whenever a knot $L$ is obtained from a knot $R$ by infection using knots $K_1,K_2,\dots$ there is a cobordism $E$ whose boundary is the disjoint union of the zero surgeries $M_L$, $-M_{R}$ and $-M_{K_1}$, $-M_{K_2}$ et cetera, as shown on the left-hand side of Figure~\ref{fig:mickey}.
\begin{figure}[htbp]
\setlength{\unitlength}{1pt}
\begin{picture}(350,200)
\put(-50,0){\includegraphics{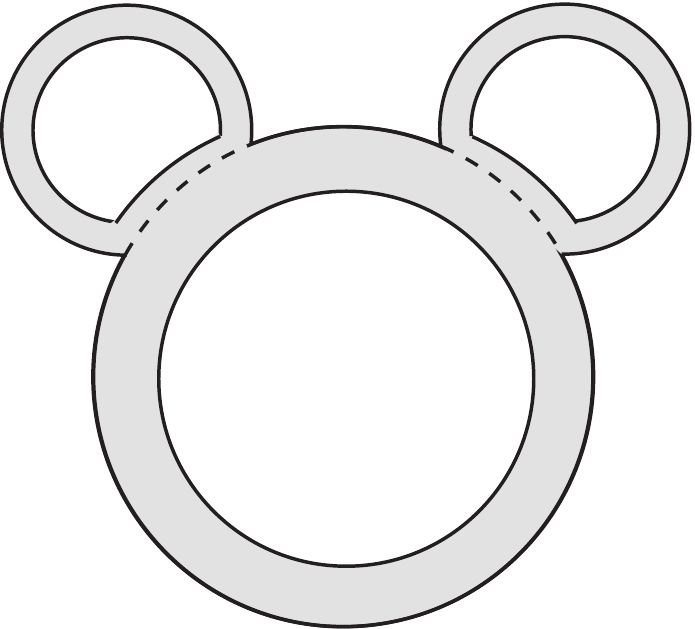}}
\put(194,0){\includegraphics{mickey.pdf}}
\put(70,100){$M_R$}
\put(88,0){$M_L$}
\put(-20,145){$M_{K_1}$}
\put(96,145){$M_{K_2}$}
\put(216,145){$M_{K^{n-1}}$}
\put(346,145){$M_{\overline{K}^{n-1}}$}
\put(327,0){$M_{K^n}$}
\put(307,93){$M_{\mathfrak{R}^{m}}$}
\put(179,73){$E_n~~\equiv$}
\put(-50,73){$E\equiv$}
\end{picture}
\caption{The cobordism}\label{fig:mickey}
\end{figure}
Therefore, since $K^n=\mathfrak{R}^{m}(K^{n-1},\overline{K^{n-1}})$, there is a cobordism $E_n$ whose boundary is the disjoint union of the zero surgeries on $K^n$, $K^{n-1}$, $\overline{K^{n-1}}$ and $\mathfrak{R}^{m}$ as shown on the right-hand side of Figure~\ref{fig:mickey} and schematically in Figure~\ref{fig:Z}. Similarly there is a cobordism $E_i$, for $1\leq i< n$ whose boundary is the disjoint union of the zero surgeries on $K^i$, $K^{i-1}$ and $R^i$. Consider $X=E_n\cup E_{n-1}\cup \overline{E}_{n-1}\cup...\cup E_1\cup \overline{E}_{1}$, gluing  $E_i$ to $E_{i-1}$ along their common boundary component $M_{K^{i-1}}$, and gluing  $\overline{E}_i$ to $\overline{E}_{i-1}$ along their common boundary component $M_{\overline{K^{i-1}}}$, as shown schematically in Figure~\ref{fig:Z}. The boundary of $X$ is a disjoint union of $M_{K^n}$, $-M_{\mathfrak{R}^{m}}$, $-M_{K^0}$, $-M_{\overline{K^0}}$ and two copies each of $\pm M_{R^{n-1}},...,\pm M_{R^1}$. For $1\leq i< n$, let $S_i$ denote the exterior of any ribbon disk in $B^4$ for the ribbon knot $R^i$. Let $S_n$  denote the exterior of any ribbon disk in $B^4$ for the ribbon knot $\mathfrak{R}^{m}$. Since Arf($K^0$)$=0$, $K^0\in  \mathcal{F}_{0}$ via some $V$ ~\cite[Section 5]{COT2}. Gluing $V$, $\overline{V}=-V$ and all the $S_i$ and $\overline{S}_i$ to $X$, we obtain a $4$-manifold, $Z$ as shown in Figure~\ref{fig:Z}. Note $\partial Z=M_{K^n}$.
\begin{figure}[htbp]
\setlength{\unitlength}{1pt}
\begin{picture}(365,283)
\put(0,15){\includegraphics{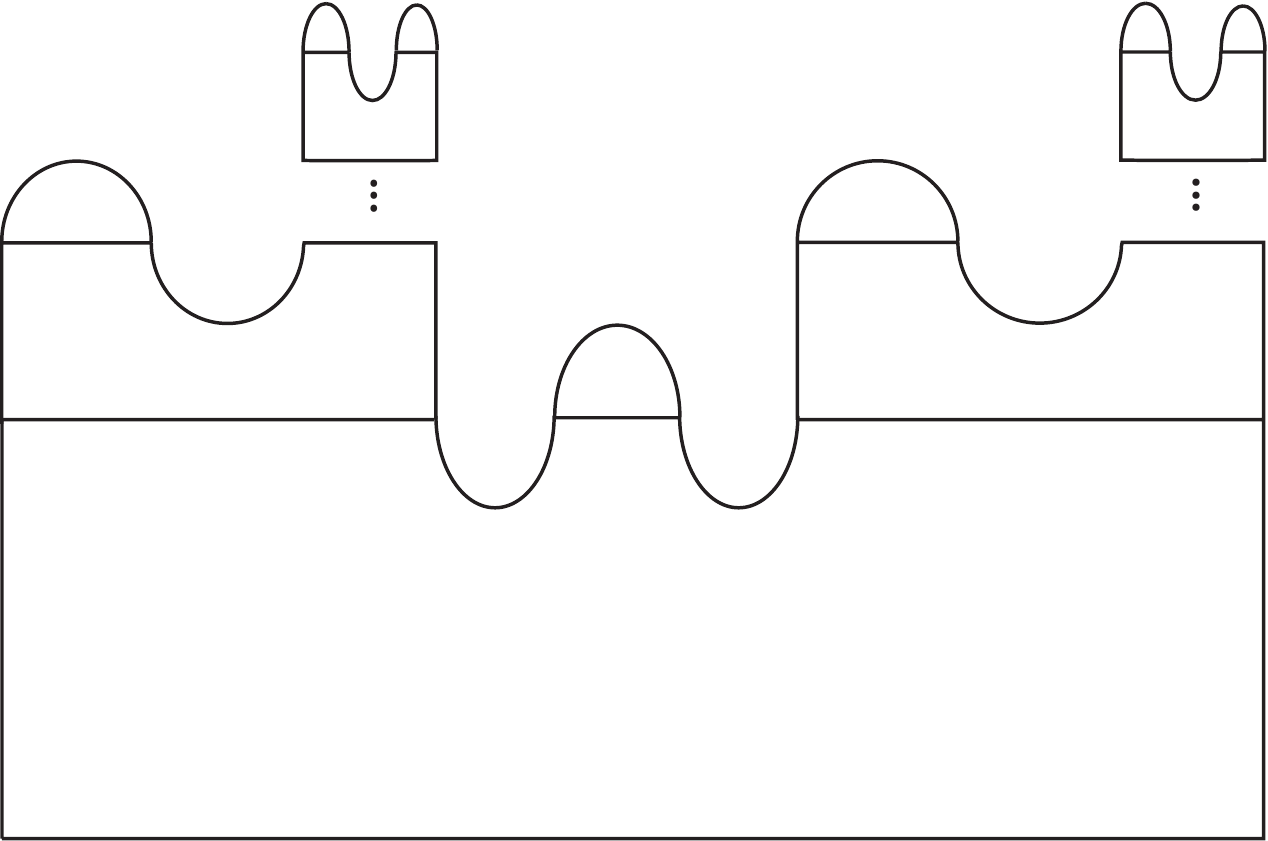}}
\put(166,127){$M_{\mathfrak{R}^{m}}$}
\put(170,146){$S_n$}
\put(167,55){$E_n$}
\put(167,4){$M_{K^n}$}
\put(241,195){$S_{n-1}$}
\put(14,193){$\overline{S}_{n-1}$}
\put(240,177){$M_{R^{n-1}}$}
\put(12,177){$M_{\overline{R}^{n-1}}$}
\put(327,177){$M_{K^{{n-2}}}$}
\put(91,177){$M_{\overline{K}^{{n-2}}}$}
\put(285,147){$E_{n-1}$}
\put(56,146){$\overline{E}_{n-1}$}
\put(326,246){$S_1$}
\put(58,247){$\overline{S}_1\to$}
\put(355,244){$V$}
\put(116,244){$\overline{V}$}
\put(340,216){$E_1$}
\put(101,216){$\overline{E}_1$}
\put(277,126){$M_{K^{n-1}}$}
\put(50,126){$M_{\overline{K}^{n-1}}$}
\end{picture}
\caption{$Z$}\label{fig:Z}
\end{figure}
We claim that,
\begin{equation}\label{eq:Knsolv}
K^n\in \mathcal{F}_{n} ~~\text{via} ~Z,
\end{equation}
and if $\mathcal{P}$ is strongly coprime to $(\Delta_m(t),q_{n-1}(t),\dots,q_1(t))$, then
\begin{equation}\label{eq:Kn+1solv}
K^n\in \mathcal{F}_{n+1}^\mathcal{P} ~~\text{via} ~Z.
\end{equation}

First, simple Mayer-Vietoris sequences together with an analysis of the homology of the $E_i$ (as given by Lemma~\ref{lem:mickeyfacts} below) imply that $H_2(Z)\cong H_2(V)\oplus H_2(\overline{V})$ since $H_2(S_i)=0$. Since $V$ is a $0$-solution, $H_2(V)$ has a basis of connected compact oriented surfaces, $\{L_j,D_j|1\leq j\leq r\}$, satisfying the conditions of Definition~\ref{def:Gnsolvable}. Similarly for $H_2(\overline{V})$.
We claim that,
\begin{equation}\label{eq:mu0}
\pi_1(V)\subset \pi_1(Z)^{(n)}
\end{equation}
and if $\mathcal{P}$ is strongly coprime to $(\Delta_m(t),q_{n-1}(t),\dots,q_1(t))$ then
\begin{equation}\label{eq:mu1}
\pi_1(V)\subset \pi_1(Z)^{(n+1)}_{\mathcal{P}}.
\end{equation}
Indeed equations \eqref{eq:mu0} and \eqref{eq:mu1} were shown inductively in the proof of \cite[Theorem 6.2, Theorem 6.5]{CHL5} using the fact that, for each $i$, the doubling operator $R^i_{\alpha_i}$ satisfies $\ell k(\alpha_i,R^i)=0$ leading to the fact that
$$
\pi_1(M_{K^{i-1}})\subset \pi_1(E_i)^{(1)}.
$$
Then,
$$
\pi_1(L_j)\subset\pi_1(V)\subset \pi_1(Z)^{(n)},
$$
and if $\mathcal{P}$ is strongly coprime to $(\Delta_m(t),q_{n-1}(t),\dots,q_1(t))$,
$$
\pi_1(L_j)\subset\pi_1(V)\subset \pi_1(Z)^{(n+1)}_{\mathcal{P}},
$$
and similarly for $\pi_1(D_j)$. The same holds for $\overline{V}$. This would complete the verification of claims ~(\ref{eq:Knsolv}) and ~(\ref{eq:Kn+1solv}) since $\{L_j,D_j\}$ (together with their counterparts in $\overline{V}$ would then satisfy the criteria of Definition~\ref{def:Gnsolvable}.

This concludes our sketch of the proof as given in \cite[Theorem 6.5]{CHL5}. We include the relevant result about the elementary topology of the cobordism $E$. We will need several of these properties in later proofs.

\begin{lem}[{\cite[Lemma 2.5]{CHL3}}]\label{lem:mickeyfacts} With regard to $E$ on the left-hand side of Figure~\ref{fig:mickey}, the inclusion maps induce
\begin{itemize}
\item [(1)] an epimorphism $\pi_1(M_L)\to \pi_1(E)$ whose kernel is the normal closure of the longitudes of the infecting knots $K_i$ viewed as curves $\ell_i\subset S^3-K_i\subset M_L$;
\item [(2)] isomorphisms $H_1(M_L)\to H_1(E)$ and $H_1(M_R)\to H_1(E)$;
\item [(3)] and isomorphisms $H_2(E)\cong H_2(M_L)\oplus_i H_2(M_{K_i})\cong H_2(M_R)\oplus_i H_2(M_{K_i})$.
\item [(4)] The meridian of $K$, $\mu_K\subset M_{K}$ is isotopic in $E$ to both $\alpha\subset M_R$ and to the longitudinal push-off of $\alpha$, often called $\alpha\subset M_L$ by abuse of notation.
\item [(5)] The longitude of $K$, $\ell_K\subset M_{K}$ is isotopic in $E$ to the reverse of the meridian of $\alpha$, $(\mu_{\alpha})^{-1}\subset M_L$ and to the longitude of $K$ in $S^3-K\subset M_L$ and to the reverse of the meridian of $\alpha$, $(\mu_{\alpha})^{-1}\subset M_R$ (the latter bounds a disk in $M_R$).
\end{itemize}
\end{lem}

\end{proof}

\section{Proof of Theorem~\ref{thm:main}}\label{sec:proof}

The proof of Theorem~\ref{thm:main} will occupy the remainder of the paper.

\newtheorem*{thm:main}{Theorem~\ref{thm:main}}
\begin{thm:main} Consider knots $K^n$, $n\geq 2$ as in Figure~\ref{fig:examplesKn}
$$
K^n\equiv \mathfrak{R}^{m}_{\eta_1,\eta_2}(K^{n-1},\overline{K^{n-1}}),
$$
where $K^{n-1}$ is the result of applying a composition of $n-1$ doubling operators, $R_{\alpha_{n-1}}^{n-1}\circ\dots\circ R_{\alpha_1}^1$ to some Arf invariant zero input knot $K^0$. Suppose additionally that
\begin{itemize}
\item [1.] $m\neq 0$;
\item [2.] for each $i$, $\alpha_i$ generates the rational Alexander module of $R^i$, and this module is non-trivial;
\item [3.] $|\rho_0(K^0)|$, the average Levine-Tristram signature of $K^0$, is greater than twice the sum of the Cheeger-Gromov constants of the ribbon knots $\mathfrak{R}^{m}$, $R^{1},\dots, R^{n-1}$ (see Section~\ref{sec:rhoinvs}).
\end{itemize}
If $\mathcal{P}$ is the $n$-tuple of Alexander polynomials of the knots $(\mathbb{E}^m, R^{n-1},\dots, R^1)$, then
$$
K^n\notin \mathcal{F}_{n.5}^{\sP}.
$$
\end{thm:main}

\begin{proof}[Proof of Theorem~\ref{thm:main}] We assume that
$$
\mathcal{P}=(p_1(t),\dots,p_n(t))=(\Delta_{m},q_{n-1}(t),\dots,q_1(t))
$$
is the $n$-tuple of Alexander polynomials of the knots $(\mathbb{E}^m, R^{n-1},\dots, R^1)$. Suppose that $K^n\in \mathcal{F}_{n.5}^{\sP}$. Let $V$ be the putative $(n.5,\mathcal{P})$-solution. We will derive a contradiction.

Let $W_0$ be the $4$-manifold (refer to Figure~\ref{fig:W_0}) obtained from $V$ by adjoining the cobordisms $E_n$, $E_{n-1}$, $\overline{E}_{n-1},\dots E_1, \overline{E}_1$ as defined in the proof of Theorem~\ref{thm:KnliesinQplusone}. For specificity, set
\begin{align*}
&W_n=V,\\ &W_{n-1}=W_n\cup E_n,\\ &W_{n-2}=W_{n-1}\cup E_{n-1}\cup\overline{E}_{n-1},\\&\vdots\\ &W_0=W_1\cup E_1\cup \overline{E}_1.
\end{align*}
Note that, unlike in the manifold $Z$ of Figure~\ref{fig:Z}, we do not cap off the zero surgeries on the various ribbon knots. Thus the boundary of $W_0$ is the disjoint union of the zero surgeries on the ribbon knots $\mathfrak{R}^{m}$, $R^{n-1}, \dots, R^{1}$, $\overline{R}^{{n-1}}, \dots, \overline{R}^{1}$, together with the zero surgeries on $K^0$, $\overline{K}^0$.
\begin{figure}[htbp]
\setlength{\unitlength}{1pt}
\begin{picture}(365,303)
\put(0,0){\includegraphics{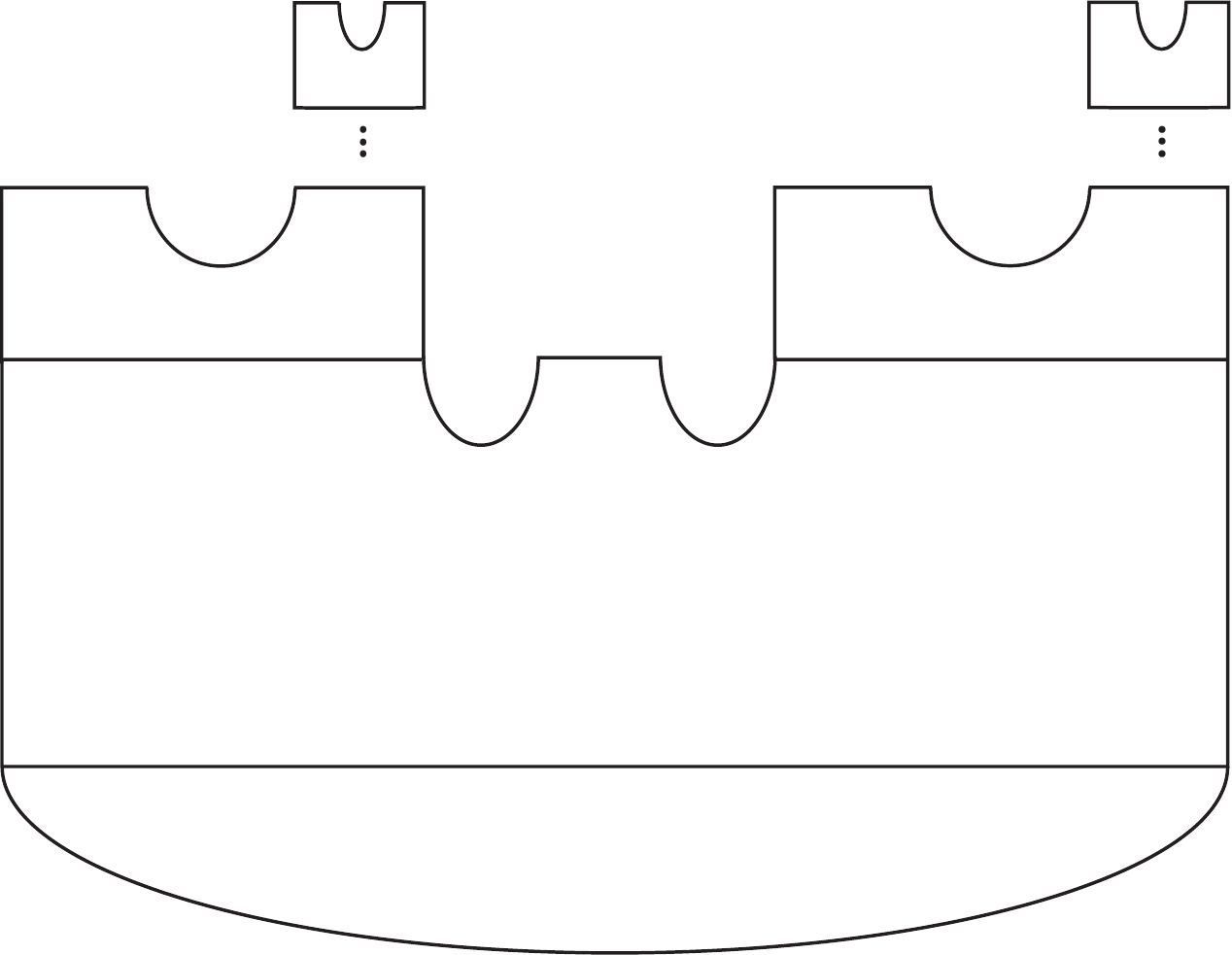}}
\put(169,183){$M_{\mathfrak{R}^{m}}$}
\put(170,100){$E_n$}
\put(172,10){$V$}
\put(323,289){$M_{R_1}$}
\put(84,290){$M_{\overline{R}_1}$}
\put(354,289){$M_{K^0}$}
\put(114,290){$M_{\overline{K}^0}$}
\put(167,45){$M_{K^n}$}
\put(240,236){$M_{R^{n-1}}$}
\put(10,236){$M_{\overline{R}^{n-1}}$}
\put(328,218){$M_{K^{{n-2}}}$}
\put(92,218){$M_{\overline{K}^{{n-2}}}$}
\put(289,190){$E_{n-1}$}
\put(58,190){$\overline{E}_{n-1}$}
\put(340,256){$E_1$}
\put(101,256){$\overline{E}_1$}
\put(278,167){$M_{K^{n-1}}$}
\put(53,167){$M_{\overline{K}^{n-1}}$}
\end{picture}
\caption{$W_{0}$}\label{fig:W_0}
\end{figure}

Below we will define a commutator series $\{\pi^{(n)}_{\mathcal{S}}\}$ that is slightly larger than the derived series localized at $\mathcal{P}$. In particular,
\begin{equation}\label{eq:PinPprime}
\pi_1(W_0)^{(n+1)}_\mathcal{P}\subset\pi_1(W_0)^{(n+1)}_{\mathcal{S}}.
\end{equation}
Then we consider the coefficient system on $W_0$ given by the projection
$$
\phi:\pi_1(W_0)\to\pi_1(W_0)/\pi_1(W_0)^{(n+1)}_{\mathcal{P}}\to\pi_1(W_0)/\pi_1(W_0)^{(n+1)}_{\mathcal{S}}.
$$
The bulk of the proof ($14$ pages!) will be to show that:
\begin{equation}\label{eq:nontrivialonK}
\text{the restriction of} ~\phi ~\text{to} ~\pi_1(M_{K^0}\subset \partial W_0) ~\text{factors non-trivially through} ~\Z; ~\text{and}
\end{equation}
\begin{equation}\label{eq:trivialonKbar}
\text{the restriction of} ~\phi ~\text{to} ~\pi_1(M_{\overline{K}^0}\subset \partial W_0)  ~\text{is zero.}
\end{equation}

We now show that \eqref{eq:PinPprime}, \eqref{eq:nontrivialonK} and \eqref{eq:trivialonKbar} imply Theorem~\ref{thm:main}.  Consider the von Neumann signature defect of ($W_0,\phi$):
$$
\sigma^{(2)}(W_0,\phi)-\sigma(W_0).
$$
By the additivity of these signatures (property ($5$) of Proposition~\ref{prop:rhoprops}), this quantity is the sum of the signature defects for $V$ and those of the $E_i$ and $\overline{E}_i$. Note that the coefficient system on $\pi_1(V)$ factors
$$
\pi_1(V)\to \frac{\pi_1(V)}{\pi_1(V)^{(n+1)}_\mathcal{P}}\to \frac{\pi_1(W_0)}{\pi_1(W_0)^{(n+1)}_\mathcal{P}}\to \frac{\pi_1(W_0)}{\pi_1(W_0)^{(n+1)}_\mathcal{S}},
$$
where we used Theorem~\ref{thm:derivedlocpfunctorial} to establish the second map and we use \eqref{eq:PinPprime} for the third map. Thus, since $V$ is an $(n.5,\mathcal{P})$-solution, the signature defect of $V$ vanishes by Theorem~\ref{thm:generalsignaturesobstruct}. All of the signature defects of the $E_i$ vanish by \cite[Lemma 2.4]{CHL3} (essentially because $H_2(E)$ comes from $H_2(\partial E)$). Therefore the signature defect vanishes for $W_0$. On the other hand, by Section~\ref{sec:rhoinvs},
$$
\sigma^{(2)}(W_0,\phi)-\sigma(W_0)=\rho(\partial W_0,\phi).
$$
Hence
$$
0=\rho(M_{\mathfrak{R}^{m}}, \phi)+\dots+\rho(M_{R^{1}},\phi)+ \rho(M_{\overline{R^{1}}},\phi) + \rho(M_{K^0},\phi)+ \rho(M_{\overline{K^0}},\phi).
$$
By \eqref{eq:nontrivialonK} and properties ($1$) and ($3$) of Proposition~\ref{prop:rhoprops},
$$
\rho(M_{K^0},\phi)=\rho_0(K^0);
$$
while by \eqref{eq:trivialonKbar} and properties ($1$) and ($2$) of Proposition~\ref{prop:rhoprops}
$$
\rho(M_{\overline{K}^0},\phi)=0.
$$
But, by choice, $|\rho_0(K^0)|$ is greater than twice the sum of the Cheeger-Gromov constants of the $3$-manifolds $M_{\mathfrak{R}^{m}},\dots, M_{R^{1}}$, which is a contradiction (see property ($6$) of Proposition~\ref{prop:rhoprops}).

Therefore the proof of Theorem~\ref{thm:main} is reduced to defining a commutator series $\{\pi^{(n)}_{\mathcal{S}}\}$ such that \eqref{eq:PinPprime}, \eqref{eq:nontrivialonK} and \eqref{eq:trivialonKbar} hold.

The commutator series $\pi^{(j)}_\mathcal{S}$ will be defined only for the groups $\pi=\pi_1(W_i)$, because we need not be concerned with any other groups. It will be defined exactly as in Definition~\ref{def:defderivedlocalp} except that the sequence of right divisor sets $S_1,...,S_n$ will be slightly different than those of Definition~\ref{def:Sdefderivedlocalp}. We now define $S_1,...,S_n$. In these definitions $\pi$ is the fundamental group of one of the $W_i$.

We define
$$
S_1=S_{1}(\pi)=S_{p_1}=S_{p_1}(\pi)=\{ q_1(\mu)...q_r(\mu)~|~(p_1(t),q_j(t))=1; ~\pi/\pi^{(1)}\cong\langle\mu\rangle\}.
$$
(Note that $\pi^{(1)}=\pi^{(1)}_{r}=\pi^{(1)}_{\mathcal{P}}=\pi^{(1)}_{\mathcal{S}}$.) Before defining the other $S_i$ we make a few remarks. Since $p_1(t)$ is a knot polynomial, $p_1(t)\doteq p_1(t^{-1})$, so $S_1$ is closed (up to units) under the natural involution. In fact, since $p_1(t)=\Delta_m(t)$ is the Alexander polynomial of $\mathbb{E}^m$, $p_1(t)$ is prime. Hence one sees that
$$
S_1=\Q[\mu,\mu^{-1}]-\langle \Delta_m(\mu)\rangle.
$$
Therefore for any $\Q[\mu^{\pm 1}]$-module $M$, $MS^{-1}_1=M_{\langle \Delta_m\rangle}$, the \textbf{classical localization of $\boldsymbol{M}$ at the prime ideal $\boldsymbol{\langle \Delta_m \rangle}$}.

Therefore, by \eqref{eqnarray:defderivedp},
\begin{equation}\label{eq:pi2S}
\pi^{(2)}_{\mathcal{S}}=\pi^{(2)}_{\mathcal{P}}\equiv \ker \left(\pi^{(1)}\to \frac{\pi^{(1)}}{[\pi^{(1)},\pi^{(1)}]}\otimes \Q[\mu,\mu^{-1}] S_{1}^{-1}\equiv \A(W)S^{-1}_{1}\equiv \A(W)_{(\Delta)}\right),
\end{equation}
where $\boldsymbol{\A(W)_{(\Delta)}}$ is the classical localization of $\A(W)$ at the prime $\langle \Delta_m \rangle$. (If $W$ is any space with $\pi_1(W)=\pi$ and $H_1(W)\cong\Z$ then by its integral Alexander module, denoted $\mathcal{A}^\Z(W)$ we mean $H_1(W;\Z[\mu,\mu^{-1}])\cong \pi^{(1)}/\pi^{(2)}$. By its rational Alexander module, denoted $\mathcal{A}(W)$, we mean $H_1(W;\Q[\mu,\mu^{-1}])$.)

Now let $\G\equiv\pi/\pi^{(2)}_{\mathcal{S}}\equiv \pi/\pi^{(2)}_{\mathcal{P}}$ and $A=\pi^{(1)}/\pi^{(2)}_{\mathcal{S}}\equiv \pi^{(1)}/\pi^{(2)}_{\mathcal{P}}\subset \G$. Thus $\G$ is the semidirect product of the abelian group $A$ with $\pi/\pi^{(1)}\cong\Z$. Note that the circle $\eta_2$ (see Figure~\ref{fig:frakR}) represents an element of $\pi_1(M_{K^n})^{(1)}$ and hence, under inclusion, an element of $\pi^{(1)}$ for each of the groups $\pi=\pi_1(W_i)$ under consideration. Hence, for any $\pi$, $\eta_2$ has an unambiguous interpretation as an element of $A$. By abuse of notation we allow $\eta_2$ to stand for its image in any of the appropriate groups. Recall that a set $S\subset\G$ is \textbf{$\boldsymbol{\G}$-invariant} if $gsg^{-1}\in S$ for all $s\in S$ and $g\in \G$. Note that the set $\{\mu^i\eta_2\mu^{-i}~|~i\in\Z\}$ is $\G$-invariant where $\mu\in\G$ generates $\pi/\pi^{(1)}$. Then we define the other $S_n$ as follows:

\begin{defn}\label{defn:S_2} Let $S_2=S_2(\pi)\subset \Q[\pi^{(1)}/\pi^{(2)}_{\mathcal{S}}]\subset \Q[\pi/\pi^{(2)}_{\mathcal{S}}]$ be the multiplicative set generated by
$$
\{~q(a)~| ~\widetilde{(q,p_2)}=1, ~q(1)\neq 0, ~a\in A\}\cup \{ ~p_2(\mu^i\eta_2\mu^{-i})~|~i\in \Z\};
$$
and for $2<i\leq n$ let
\begin{eqnarray}\label{eq:higherSn}
S_{n}=S_{n}(\pi)=\{ q_1(a_1)...q_r(a_r)~|~\widetilde{(p_n,q_j)}=1; ~q_j(1)\neq 0; ~a_j\in \pi^{(n-1)}_\mathcal{S}/\pi^{(n)}_\mathcal{S}\}.
\end{eqnarray}
\end{defn}

Since $S_2$ is a multiplicative subset of $\Q A$ that is $\G$-invariant, it is a right divisor set of $\Q\G$ by \cite[Proposition 4.1]{CHL5}. Therefore Definition~\ref{def:defderivedlocalp} applies to give a  partially defined commutator series $\{\pi^{(i)}_{\mathcal{S}}\}$. Since $p_2(t)=q_{n-1}(t)$ is a knot polynomial, $p_2(t)\doteq p_2(t^{-1})$. Thus $S_2$ is closed (up to units) under the natural involution.

\begin{lem}\label{lem:pipinpis} For each $\pi$ and each $0\leq i\leq n+1$
\begin{equation}\label{eq:PinPprime2}
\pi^{(i)}_{\mathcal{P}}\subset\pi^{(i)}_{\mathcal{S}}.
\end{equation}
\end{lem}
\begin{proof} The proof is by induction on $i$. By Definition~\ref{def:defderivedlocalp}, $\pi^{(1)}_\mathcal{P}=\pi^{(1)}_r=\pi^{(1)}_\mathcal{S}$, so the Lemma is true for $i=0,1$. Suppose it is true for all values up to some fixed $i\geq 1$. Let $j:\pi\to \pi$ be the identity map. By ~\cite[Proposition 3.2]{CHL5}, it suffices to show that the induced ring map
$$
j_*:\Z[\pi/\pi^{(i)}_\mathcal{P}]\to \Z[\pi/\pi^{(i)}_\mathcal{S}]
$$
has the property that $j_*(S_{p_i}(\pi))\subset S_i(\pi)$. For $i=1$, $j_*$ is the identity map and $S_1(\pi)$ is, by definition, identical to $S_{p_1}(\pi)$. It follows that $\pi^{(2)}_\mathcal{P}=\pi^{(2)}_\mathcal{S}$ as already observed in ~\eqref{eq:pi2S}. Thus, for $i=2$, $j_*$ is again the identity map and, by Definitions~\ref{defn:S_2} and ~\ref{def:Sdefderivedlocalp}, $S_2(\pi)$ strictly contains $S_{p_2}(\pi)$. For $i>2$, the map $j_*$, although induced by the identity, will be a surjection with non-zero kernel. Nonetheless, by the inductive hypothesis, $j$ induces a homomorphism
$$
j_*:\pi^{(i-1)}_\mathcal{P}/\pi^{(i)}_\mathcal{P}\to \pi^{(i-1)}_\mathcal{S}/\pi^{(i)}_\mathcal{S}.
$$
Recall from Definition~\ref{def:Sdefderivedlocalp} that
$$
S_{p_i}(\pi)=\{ q_1(a_1)...q_r(a_r)~|~\widetilde{(p_i,q_j)}=1; ~q_j(1)\neq 0; ~a_j\in \pi^{(i-1)}_\mathcal{P}/\pi^{(i)}_\mathcal{P}\},
$$
which is the multiplicative set generated by the described set of polynomials $q(a)$. If $q(a)$ is any such polynomial then $j_*(q(a))=q(j_*(a))$ and since
$$
a\in \frac{\pi^{(i-1)}_\mathcal{P}}{\pi^{(i)}_\mathcal{P}},~~ j_*(a)\in \frac{\pi^{(i-1)}_\mathcal{S}}{\pi^{(i)}_\mathcal{S}}.
$$
Thus, upon examining ~\eqref{eq:higherSn}, we see that $q(j_*(a))\in S_i(\pi)$. Hence $j_*(S_{p_i}(\pi))\subset S_i(\pi)$ as desired.
\end{proof}

In particular this establishes \eqref{eq:PinPprime}.

\begin{lem}\label{lem:Sisfunctorial} The commutator series $\{\pi^{(i)}_\mathcal{S}\}$ is functorial with respect to any inclusion, $W_i\to W_j$, where $i>j$.
\end{lem}

\begin{proof} Note that any such inclusion induces an isomorphism $H_1(W_i)\cong H_1(W_j)\cong\Z=\langle \mu\rangle$.  If $\pi^{(i)}_\mathcal{S}$ were actually the polarized derived series localized at $\mathcal{P}$, then the functoriality would follow directly from  our Theorem~\ref{thm:derivedlocpfunctorial} \cite[Thm. 4.16]{CHL5}. But since $\pi^{(i)}_\mathcal{S}$ is slightly different, we must actually repeat some of the proof of \cite[Thm. 4.16]{CHL5}. Suppose $A=\pi_1(W_i), B=\pi_1(W_j)$ and $\psi:A\to B$ is induced by inclusion. We show, by induction on $i$, that $\psi(A^{(i)}_\mathcal{S})\subset B^{(i)}_\mathcal{S}$. This holds for $i=0$ so suppose it holds for $i=n$. We will show that $\psi(A^{(n+1)}_\mathcal{S})\subset B^{(n+1)}_\mathcal{S}$. The induction hypothesis guarantees that, for each $1\leq k\leq n$, $\psi$ induces a homomorphism of pairs
$$
\psi:(A/A^{(k)}_\mathcal{S}, A^{(k-1)}_\mathcal{S}/A^{(k)}_\mathcal{S})\to (B/B^{(k)}_\mathcal{S},B^{(k-1)}_\mathcal{S}/B^{(k)}_\mathcal{S}).
$$
By \cite[Prop.3.2]{CHL5} (or by examining ~\eqref{eqnarray:defderivedp}) it suffices to show that this map satisfies
\begin{equation}\label{eq:functor}
\psi(S_k(A))\subset S_k(B)
\end{equation}
for each $1\leq k\leq n$. First consider $k=1$. Recall that
$$
S_{1}(A)=\{ q_1(\mu)...q_r(\mu)~|~(p_1(t),q_j(t))=1; ~A/A^{(1)}\cong\langle\mu\rangle\}\subset \Q[A/A^{(1)}].
$$
Since $\psi$ induces an isomorphism $\psi:A/A^{(1)}\to B/B^{(1)}$, $\psi(\mu)=\pm\mu$. By choosing generators once and for all, we may assume that $\psi(\mu)=\mu$. So, for any such $q_j(t)$,
$$
\psi(q_1(\mu)\dots q_r(\mu))=q_1(\psi(\mu))\dots q_r(\psi(\mu))=q_1(\mu)\dots q_r(\mu)\in S_{1}(B).
$$
This verifies ~\eqref{eq:functor} for $k=1$.

Now suppose $k>1$. Recall that
$$
S_{k}(A)=\{ q_1(a_1)...q_r(a_r)~|~\widetilde{(p_n,q_j)}=1; ~q_j(1)\neq 0; ~a_j\in A^{(k-1)}_\mathcal{S}/A^{(k)}_\mathcal{S}\}.
$$
So, for any such $q_j(t)$,
$$
\psi(q_1(a_1)\dots q_r(a_r))=q_1(\psi(a_1))\dots q_r(\psi(a_r))\in S_{k}(B),
$$
since $\psi(a_j)\in B^{(k-1)}_\mathcal{S}/B^{(k)}_\mathcal{S}$.

Thus $\psi(S_{k}(A))\subset S_{k}(B)$.
\end{proof}

\begin{subsection}{Establishing \eqref{eq:nontrivialonK} and \eqref{eq:trivialonKbar}}\label{subsec:crucialsteps}\

Since $\pi_1(M_{K^0})\subset \pi_1(W_0)$ is normally generated by its meridian, $\mu_0$, and $\pi_1(M_{\overline{K}^0})$ is normally generated by its meridian, (that we denote) $\overline{\mu}_0$, the case $i=0$ of the following Proposition will establish \eqref{eq:nontrivialonK} and \eqref{eq:trivialonKbar}. Therefore the rest of the paper will be spent establishing Proposition~\ref{prop:inductivestep}.

\begin{prop}\label{prop:inductivestep} For any $i$, $0\leq i\leq n-2$, $\mu_{i}=\alpha_{i+1}$ is non-trivial, while $\overline{\mu}_{i}=\overline{\alpha}_{i+1}$ is trivial in
$$
\frac{\pi_1(W_{i})^{(n-i)}}{\pi_1(W_{i})^{(n-i+1)}_{\mathcal{S}}}.
$$
\end{prop}

To clarify the notation of this proposition, recall that, for $0\leq i\leq n-1$, $\partial W_{i}$ contains the disjoint union of the zero surgeries on the knots $K^{i}$ (refer to the schematic Figure~\ref{fig:Eidashed}), and $\overline{K}^{i}$. Let $\mu_{i}$ and $\overline{\mu}_{i}$ denote the meridians of $K^{i}$ and $\overline{K}^{i}$ in these copies of $M_{K^{i}}$ and $M_{\overline{K}^{i}}$ respectively.  Also recall that $K^{i+1}=R^{i+1}_{\alpha_{i+1}}(K^{i})$ for some circle $\alpha_{i+1}$ that generates the Alexander module of $R^{i+1}$; and $\overline{K}^{i+1}=\overline{R}^{i+1}_{\alpha_{i+1}}(\overline{K}^{i})$. Let $\alpha_{i+1}$ denote (a push-off of) this circle in $M_{K^{i+1}}\subset \partial W_{i+1}$ (referring  to Figure~\ref{fig:Eidashed}); and let $\overline{\alpha}_{i+1}$ denote (a push-off of) the other copy of $\alpha_{i+1}$ in $M_{\overline{K}^{i+1}}\subset \partial W_{i+1}$. Note that, by property ($4$) of Lemma~\ref{lem:mickeyfacts}, $\mu_{i}$ is isotopic to $\alpha_{i+1}$ in $E_{i+1}$ and $\overline{\mu}_{i}$ is isotopic to $\overline{\alpha}_{i+1}$ in $\overline{E}_{i+1}$. Hence $\mu_{i}=\alpha_{i+1}$ and $\overline{\mu}_{i}=\overline{\alpha}_{i+1}$ as elements of $\pi_1(W_{i})$.

\begin{figure}[htbp]
\setlength{\unitlength}{1pt}
\begin{picture}(218,217)
\put(0,0){\includegraphics{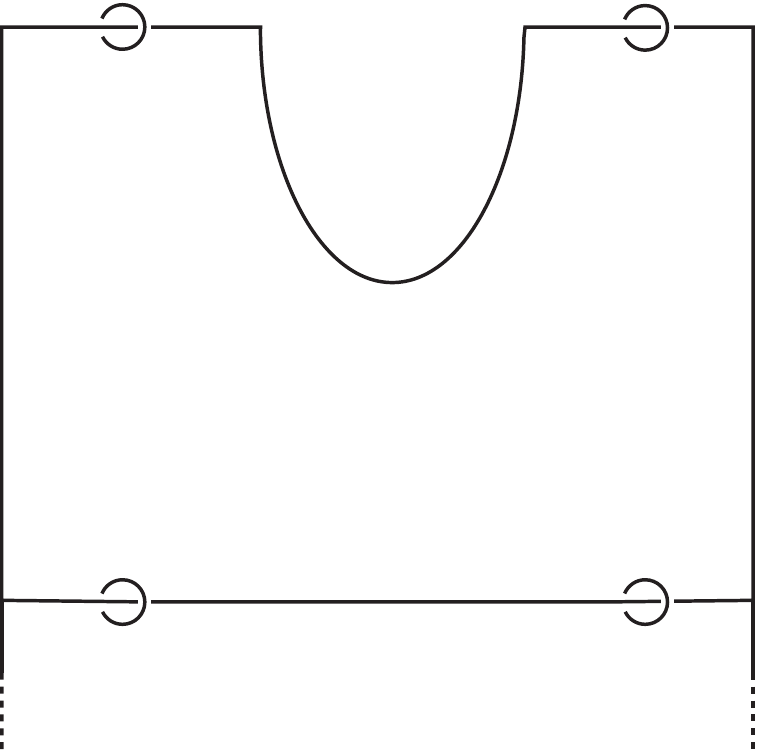}}
\put(179,194){$\mu_i$}
\put(25,194){$\alpha_{i+1}$}
\put(100,100){$E_{i+1}$}
\put(105,6){$W_{i+1}$}
\put(105,50){$M_{K^{i+1}}$}
\put(172,56){$\alpha_{i+1}$}
\put(23,56){$\mu_{i+1}$}
\put(229,209){$M_{K^i}$}
\put(-35,209){$M_{R^{i+1}}$}
\end{picture}
\caption{$W_{i}$}\label{fig:Eidashed}
\end{figure}

\begin{proof}[Proof of Proposition~\ref{prop:inductivestep}] The proof is by reverse induction on $i$, starting with $i=n-2$.

Before proving the base case $i=n-2$, we need to work out the ``pre-base-case'', $i=n-1$, where the situation is slightly different. Note that $\alpha_n$ and $\overline{\alpha}_n$ are what we have previously called $\eta_1$ and $\eta_2$ respectively.

\begin{lem}\label{lem:easystep} $\mu_{n-1}=\eta_1$ and $\overline{\mu}_{n-1}=\eta_2$ are both non-trivial in
$$
\frac{\pi_1(W_{n-1})^{(1)}}{\pi_1(W_{n-1})^{(2)}_{\mathcal{S}}}.
$$
\end{lem}

\begin{proof} Throughout the proof of this lemma we abbreviate $W=W_{n-1}$, $\pi=\pi_1(W_{n-1})$ and $\Delta=\Delta_{n}$.  We make use of the fact that the integral and rational Alexander modules of a knot agree with those of its zero-framed surgery. Specifically we use $\mathcal{A}(K)$ to denote both the rational Alexander module of $K$ and that of $M_K$. The inclusion maps induce a commutative diagram of maps between integral and rational Alexander modules as shown:
$$
\begin{diagram}\label{diag:Alexandermods}
\node{\A^\Z (\mathbb{E}^m)}\arrow{e,t}{i}\arrow{s,r}{i'_1}\node{\A^\Z (\mathfrak{R}^{m})}\arrow{e,t}{j_*}\arrow[2]{s,r}{i_2}\node{\A^\Z (V)}\arrow{e,t}{k_*}\arrow[2]{s,r}{i_3}
\node{\A^\Z (W)}\arrow{e,t}{}\arrow[2]{s,r}{i_4}\node{\frac{\pi^{(1)}}{\pi_{\mathcal{S}}^{(2)}}}\arrow{ssw,r}{i_5}\\
\node{\A (\mathbb{E}^m)}\arrow{s,r}{i''_{1}}\\
\node{\A (\mathbb{E}^m)_{(\Delta)}}\arrow{e,t}{i}\node{\A (\mathfrak{R}^{m})_{(\Delta)}}\arrow{e,t}{j_*}\node{\A (V)_{(\Delta)}}\arrow{e,t}{k_*}\node{\A (W)_{(\Delta)}}
\end{diagram}
$$
Notice that $\A^\Z (\mathfrak{R}^{m})\cong\A^\Z (K^{n})\cong\A^\Z (M_{K^{n}})\cong\A^\Z (\partial V)$.
The maps $j_*$ and $k_*$ are induced by inclusion. The map $i$ is induced by the connected sum decomposition, where here $\mathbb{E}^m$ denotes the ``left-hand'' copy in $\mathfrak{R}^{m}\equiv \mathbb{E}^m~\#~\mathbb{E}^m$. The existence and injectivity of $i_5$ is given by \eqref{eq:pi2S}. Since the $\eta_i$ represent elements in the Alexander module of $\mathbb{E}^m$, it suffices to show that the composition in the top row is injective. For this it suffices to show that the composition $k_*\circ j_*\circ i\circ i''_{1}\circ i'_{1}$ is a monomorphism. Since it is well known that the integral Alexander modules $\A^\Z (\mathbb{E}^m)\cong \A^\Z (S^3-\mathbb{E}^m)$ are $\Z$-torsion-free, $i'_1$ is injective. Since $\A (\mathbb{E}^m)$ is a $\Delta$-torsion module, $i''_1$ is injective. Under the connected sum decomposition the localized Alexander module of $\mathfrak{R}^{m}$ decomposes as the direct sum of the localized Alexander modules of its summands $\mathbb{E}^m$. The Blanchfield form decomposes similarly. Hence $i$ is injective. Now consider the map $j_*$ induced by the inclusion $\partial V\hookrightarrow V$.
$$
j_*:\A (\partial V)_{(\Delta)}\cong\A (\mathfrak{R}^{m})_{(\Delta)}\equiv H_1(M_{\mathfrak{R}^{m}};\Q[t,t^{-1}]S_{p_1}^{-1})\to H_1(V;\Q[t,t^{-1}]S_{p_1}^{-1})\equiv \A (V)_{(\Delta)}.
$$
Since $V$ is an $(n.5,\mathcal{P})$-solution for $\partial V$, and $\pi^{(i)}_{\mathcal{P}}\subset\pi^{(i)}_{\mathcal{S}}$, $V$ is an $(n.5,\mathcal{S})$ solution, so it is certainly a $(1,\mathcal{S})$-solution. Consider the coefficient system $\psi:\pi_1(V)\to \pi_1(V)/\pi_1(V)^{(1)}_{\mathcal{S}}\cong \Z$ (recall $G^{(1)}_S=G^{(1)}_r$ for any group $G$). Then ~\cite[Theorem 7.15]{CHL3} applies to say that the kernel of $j_*$ is isotropic with respect to the classical Blanchfield form on $\A (\mathfrak{R}^{m})_{(\Delta)}$. Hence the kernel, $P$, of $i\circ j_*$ is isotropic with respect to the classical Blanchfield form on $\A (\mathbb{E}^m)_{(\Delta)}$. But, since the Alexander polynomial of $\mathbb{E}^m$ is irreducible by Proposition~\ref{prop:alexpolyscoprime}, the rational Alexander module of $\mathbb{E}^m$ has no proper submodules. The case $P=\A (\mathbb{E}^m)_{(\Delta)}$ is not possible since the localized classical Blanchfield form is non-singular and $\A (\mathbb{E}^m)_{(\Delta)}\neq 0$. Thus $P=0$ so $i\circ j_*$  is injective.

It only remains to show that the lower map $k_*$ is injective (actually an isomorphism). Since localization is an exact functor, this is equivalent to showing that the inclusion map induces an isomorphism between the rational Alexander modules of $V$ and $W$. Recall that $W=W_n=V\cup E_n$. Recall from property ($1$) of Lemma~\ref{lem:mickeyfacts} applied to $E_n$, that the kernel on $\pi_1$ of the inclusion $M_{K^n}=\partial V\to E_n$ is normally generated by the longitudes of the infecting knots $K^{n-1}$ and $\overline{K}^{n-1}$ as curves in $\pi_1(M_{K^n})$. These lie in the second derived subgroups of $\pi_1(S^3-K^{n-1})$ and $\pi_1(S^3-\overline{K}^{n-1})$ respectively and so lie in the third derived subgroup of $\pi_1(M_{K^n})$ (refer to Figure~\ref{fig:infectionannotated}). Since the rational Alexander module of any space $X$ with $H_1(X)\cong \Z$ may be described as $G^{(1)}/G^{(2)}\otimes \Q$ where $G=\pi_1(X)$, this shows that the rational Alexander modules of $V$ and $W$ are isomorphic.
\end{proof}

The crucial base case, $i=n-2$, in the (reverse) inductive proof of Proposition~\ref{prop:inductivestep} is:

\begin{lem}\label{lem:crucialstep} $\mu_{n-2}=\alpha_{n-1}$ is non-trivial, while $\overline{\mu}_{n-2}=\overline{\alpha}_{n-1}$ is trivial in
$$
\frac{\pi_1(W_{n-2})^{(2)}}{\pi_1(W_{n-2})^{(3)}_{\mathcal{S}}}.
$$
\end{lem}

\begin{proof} It might be helpful to refer to Figure~\ref{fig:Eidashed} with $i=n-2$. By property ($1$) of Lemma~\ref{lem:mickeyfacts}, the kernel of the map
$$
\pi_1(W_{n-1})\to \pi_1(W_{n-1}\cup E_{n-1}\cup \overline{E}_{n-1})= \pi_1(W_{n-2})
$$
is normally generated by the longitudes, $\ell_{n-2}, \overline{\ell}_{n-2}$, of the infecting knots $K^{n-2}$ and $\overline{K}^{n-2}$ viewed as curves in $S^3\setminus K^{n-2}\subset M_{K^{n-1}}\subset \partial W_{n-1}$ and $S^3\setminus \overline{K}^{n-2}\subset M_{\overline{K}^{n-1}}\subset \partial W_{n-1}$. But of course these lie in the second derived subgroups of $\pi_1(S^3\setminus K^{n-2})$ and $\pi_1(S^3\setminus \overline{K}^{n-2})$ respectively, and so lie in the second derived subgroups of $\pi_1(M_{K^{n-1}})$ and $\pi_1(M_{\overline{K}^{n-1}})$ respectively. But, as observed in Lemma~\ref{lem:easystep}
\begin{equation}\label{eq:munminusone}
\pi_1(M_{K^{n-1}})=\langle \mu_{n-1}\rangle\subset \pi_1(W_{n-1})^{(1)},
\end{equation}
and similarly for $\pi_1(M_{\overline{K}^{n-1}})$. It follows that both $\ell_{n-2}$ and $\overline{\ell}_{n-2}$ lie the third derived subgroup of $\pi_1(W_{n-1})$ and hence lie in $\pi_1(W_{n-1})^{(3)}_{\mathcal{S}}$. Thus the inclusion $W_{n-1}\to W_{n-2}$ induces an isomorphism
$$
\frac{\pi_1(W_{n-1})}{\pi_1(W_{n-1})^{(3)}_{\mathcal{S}}}\cong \frac{\pi_1(W_{n-2})}{\pi_1(W_{n-2})^{(3)}_{\mathcal{S}}},
$$
by weak functoriality and by \cite[Prop.4.7]{CHL5}.

Therefore, to prove Lemma~\ref{lem:crucialstep}, it suffices to let $\pi=\pi_1(W_{n-1})$, and show that $\alpha_{n-1}$ is non-trivial in $\pi^{(2)}/\pi^{(3)}_{\mathcal{S}}$ and that $\overline{\alpha}_{n-1}$ is trivial in $\pi^{(2)}/\pi^{(3)}_{\mathcal{S}}$. Throughout the rest of the proof of Lemma~\ref{lem:crucialstep}, we will abbreviate $W=W_{n-1}$, $\pi=\pi_1(W_{n-1})$, $J=K^{n-1}$ and $\overline{J}=\overline{K}^{n-1}$. Thus $\partial W=M_{\mathfrak{R}^{m}}\cup M_J \cup M_{\overline{J}}$.

Consider the following commutative diagram (which we justify below) where $\G=\pi/\pi^{(2)}_\mathcal{S}$ and $\mathcal{R}=\Q\G S_2^{-1}$. Since we may view $\alpha_{n-1}\in\pi_1(M_J)^{(1)}$ and $\overline{\alpha}_{n-1}\in \pi_1(M_{\overline{J}})^{(1)}$, we have reduced Lemma~\ref{lem:crucialstep} to showing that $\alpha_{n-1}$ is \emph{not} in the kernel of the top row of the diagram while $\overline{\alpha}_{n-1}$ \emph{does} lie in this kernel.
$$
\begin{diagram}\label{diagram-pirelatehomology}\dgARROWLENGTH=.8em
\node{\pi_1(M_J)^{(1)}\oplus\pi_1(M_{\overline{J}})^{(1)}}\arrow[2]{e,t}{j_*}\arrow{s,r}{\pi}\node[2]{\pi^{(2)}}\arrow{e,t}{\phi}\arrow{s,r}{} \node{\frac{\pi^{(2)}_\mathcal{S}}{\pi^{(3)}_\mathcal{S}}}\arrow{s,r}{j} \\
  \node{(\mathcal{A}(J)\oplus \mathcal{A}(\overline{J}))\otimes \mathcal{R}}\arrow{e,t}{\cong }\node{ H_1(M_J \cup M_{\overline{J}};\mathcal{R})} \arrow{e,t}{j_*}\node{H_1(W;\R)}\arrow{e,t}
{\cong}\node{\frac{\pi^{(2)}_\mathcal{S}}{[\pi^{(2)}_\mathcal{S},\pi^{(2)}_\mathcal{S}]}\otimes \mathcal{R}}
\end{diagram}
$$
The $j_*$ in the upper row of the diagram is justified by our observation \eqref{eq:munminusone}, which says that $\pi_1(M_J)\subset \pi^{(1)}$ and $\pi_1(M_{\overline{J}})\subset \pi^{(1)}$. Now we consider the first map in the bottom row. By Lemma~\ref{lem:easystep} the coefficient system $\pi\to \G$, when restricted to $\pi_1(M_J)$ is non-trivial:
$$
\pi_1(M_J)=\langle \mu_{n-1}\rangle \hookrightarrow \frac{\pi^{(1)}}{\pi^{(2)}_\mathcal{S}}\hookrightarrow \frac{\pi}{\pi^{(2)}_\mathcal{S}}\equiv \G,
$$
but also factors through $\pi_1(M_J)/\pi_1(M_J)^{(1)}\cong\Z$ using \eqref{eq:munminusone}. It follows that
$$
H_1(M_J;\Q\G)\cong H_1(M_J;\Q[t,t^{-1}])\otimes \Q\G\equiv \mathcal{A}(J)\otimes_{\Q[t,t^{-1}]} \Q\G,
$$
where $\Q[t,t^{-1}]$ acts on $\Q\G$ by $t\to \mu_{n-1}$ (equivalently $t\to \eta_1$). Hence
$$
H_1(M_J;\mathcal{R})\cong \mathcal{A}(J)\otimes \mathcal{R};
$$
and similarly for $\overline{J}$, where $t$ acts by $\overline{\mu}_{n-1}=\eta_2$. This explains the first map in the lower row of the diagram. To justify the last map in the lower row, recall that $H_1(W;\Z\G)$ has an interpretation as the first homology module of the $\G$-covering space of $W$. The fundamental group of this covering space is the kernel of $\pi\to \G$. Hence
$$
H_1(W;\Z\G)\cong \frac{\pi^{(2)}_\mathcal{S}}{[\pi^{(2)}_\mathcal{S},\pi^{(2)}_\mathcal{S}]}
$$
Since the Ore localization $\mathcal{R}$ is a flat $\Z\G$-module, the $\cong$ is justified. This completes the explanation of the diagram. Since, by Definitions~\ref{def:Sdefderivedlocalp} and \ref{defn:S_2},
$$
\pi^{(3)}_{\mathcal{S}}=\ker \left(\pi^{(2)}_\mathcal{S}\to \frac{\pi^{(2)}_\mathcal{S}}{[\pi^{(2)}_\mathcal{S},\pi^{(2)}_\mathcal{S}]}\to \frac{\pi^{(2)}_\mathcal{S}}{[\pi^{(2)}_\mathcal{S},\pi^{(2)}_\mathcal{S}]}\otimes \Q\G S_2^{-1}\right),
$$
it follows that the vertical map $j$ (in the diagram) is injective. Hence, to establish Lemma~\ref{lem:crucialstep}, it suffices to show that the class represented by $\alpha_{n-1}\otimes 1$ is \emph{not} in the kernel of the bottom row of the diagram while that represented by $\overline{\alpha}_{n-1}\otimes 1$ \emph{does} lie in this kernel.

Recall that $\overline{J}\equiv \overline{K}^{n-1}\equiv \overline{R}^{n-1}_{\overline{\alpha}_{n-1}}(\overline{K}^{n-2})$  where $\overline{\alpha}_{n-1}$ generates $\mathcal{A}(\overline{R}^{n-1})$ (note this implies the latter module is cyclic). Therefore $\mathcal{A}(\overline{J})\cong \mathcal{A}(\overline{R}^{n-1})$. By hypothesis, the Alexander polynomial of $R^{n-1}$ is $q_{n-1}(t)=p_2(t)$. Thus
$$
\langle \overline{\alpha}_{n-1}\rangle\cong\mathcal{A}(\overline{J})\cong\frac{\Q[t,t^{-1}]}{p_2(t)\Q[t,t^{-1}]}
$$
and
$$
\langle \overline{\alpha}_{n-1}\otimes 1\rangle\cong\mathcal{A}(\overline{J})\otimes\mathcal{R}\cong \left (\frac{\Q\G}{p_2(\eta_2)\Q\G}\right )S_2^{-1}\cong0,
$$
where the last equality holds since $p_2(\eta_2)\in S_2$, by Definition~\ref{defn:S_2} (see \cite[Thm. 4.12]{CHL5} for more detail). Therefore $\overline{\alpha}_{n-1}\otimes 1$ lies in the kernel of the bottom row of the diagram.

Suppose that $\alpha_{n-1}\otimes 1$ \emph{were} in the kernel of the bottom row of the diagram. We shall reach a contradiction. Recall that $W_{n-1}\equiv V\cup E_n$. Recall that $V$ is an $(n.5,\mathcal{P})$-solution. Since $n\geq 2$, $V$ is a $(2,\mathcal{P})$-solution. One easily checks that
$$
\frac{H_2(W_{n-1})}{i_*(H_2(\partial W_{n-1})}\cong H_2(V).
$$
Hence this group has a basis consisting of surfaces that satisfy parts ($2$) and ($3$) of Definition~\ref{def:Gnsolvable} (with $n=2$). But $W_{n-1}$ fails to satisfy part ($1$) of that definition and  $\partial W_{n-1}$ is disconnected. Such a manifold was named a \textbf{$(2,\mathcal{P})$-bordism} in \cite[Definition 7.11]{CHL5}. By ~\cite[Thms. 7.14, 7.15]{CHL5}, if $P$ is the kernel of the map
$$
j^*:H_1(M_J;\mathcal{R})\to H_1(W;\mathcal{R}),
 $$
as in the bottom row of the diagram, then $P$ is an isotropic submodule for the Blanchfield linking form on $H_1(M_J;\mathcal{R})$.  Since we have supposed that $\alpha_{n-1}\otimes 1\in P$ and since this element is a generator of $H_1(M_J;\mathcal{R})$, it would follow that this Blanchfield form were identically zero on $H_1(M_J;\mathcal{R})$. But by \cite[Lemma 7.16]{CHL5} this form is non-singular. This would imply that $H_1(M_J;\mathcal{R})$ were the zero module. This is a contradiction once we show that
\begin{equation}\label{eq:nontrivialmodule}
\mathcal{A}({J})\otimes\mathcal{R}\cong\left (\frac{\Q\G}{p_2(\eta_1)\Q\G}\right )S_2^{-1}\neq 0.
\end{equation}
This is a non-trivial result since we are dealing with a noncommutative localization.

Note that, by the hypotheses of Theorem~\ref{thm:main}, $p_2(t)=q_{n-1}(t)$ is not a unit in $\Q[t,t^{-1}]$. The map $\Z\to \G$ given by $t\to \eta_1$ is not zero by Lemma~\ref{lem:easystep}. Since $\G$ is PTFA, it is torsion-free, so $\langle\eta_1\rangle\subset\G$. Hence $\Q\G$ is a free left $\Q[\eta_1,\eta_1^{-1}]$-module on the right cosets of $\langle\eta_1\rangle\in \G$ ~\cite[Chapter 1, Lemma 1.3]{P}. Thus, upon fixing a set of coset representatives, any $x\in \Q\G$ has a unique decomposition
$$
x=\Sigma_\gamma x_\gamma\gamma ,
$$
where $x_\gamma\in\mathbb{Q}[\eta_1,\eta_1^{-1}]$ and the sum is over a set of coset representatives $\{\gamma\in \G\}$. It follows that $p_2(\eta_1)$ has no right inverse in $\Q\G$ since if $p_2(\eta_1)x=1$ then
$$
p_2(\eta_1)x=p_2(\eta_1)\Sigma_\gamma x_\gamma\gamma= \Sigma_\gamma p_2(\eta_1)x_\gamma\gamma =1.
$$
Looking at the coset $\gamma=e$ , we have $p_2(\eta_1)x_e=1$ in $\mathbb{Q}[\eta_1,\eta_1^{-1}]$, contradicting the fact that $p_2(t)$ is not a unit in $\Q[t,t^{-1}]$. Therefore, since $\Q\G$ is a domain,
$$
\frac{\Q\G}{p_2(\eta_1)\Q\G}\ncong 0.
$$

\noindent Continuing, by \cite[Corollary 3.3, p. 57]{Ste}, the kernel of
$$
\frac{\Q\G}{p_2(\eta_1)\Q\G}\to \left (\frac{\Q\G}{p_2(\eta_1)\Q\G}\right )S_2^{-1}
$$
is precisely the $S_2$-torsion submodule. Hence to establish ~\eqref{eq:nontrivialmodule}, it suffices to show that the generator of $\Q\G/p_2(\eta_1)\Q\G$ is not $S_2$-torsion. Suppose $[1]$ were $S_2$-torsion. We will show that $[1]=0$, implying that $\Q\G/p_2(\eta_1)\Q\G$ is $S_2$-torsion-free. If $[1]$ were $S_2$-torsion then $1s=p_2(\eta_1)y$ for some $s\in S_2$ and for some $y\in \Q\G$. We examine this equation in $\Q\G$.

Recall that $\G=\pi/\pi^{(2)}_\mathcal{S}$. Let $A= \pi^{(1)}/\pi^{(2)}_\mathcal{S}\lhd \G$. Since $A\subset \G$, $\mathbb{Q}\G$, viewed as a left $\mathbb{Q}A$-module, is free on the right cosets of $A$ in $\G$. Thus any $y\in \Q\G$ has a unique decomposition
$$
y=\Sigma_\gamma y_\gamma\gamma ,
$$
where the sum is over a set of coset representatives $\{\gamma\in \G\}$  and $y_\gamma\in\mathbb{Q}A$.  Therefore we have
$$
 s=p_2(\eta_1)\Sigma_\gamma y_\gamma\gamma.
$$
Recall from Definition~\ref{defn:S_2} that $s\in S_2\subset \Q A$. It follows that for each coset representative $\gamma\neq e$ we have $0=p_2(\eta_1)y_\gamma$ so $y_\gamma=0$ (note that $p_2(\eta_1)\neq 0$ since $\Q[\eta_1^{\pm 1}]\subset \Q\G$). Hence $y\in \Q A$ and we have
\begin{equation}\label{eq:storsion}
s=p_2(\eta_1)y
\end{equation}
as an equation in $\Q A$. Recall from Definition~\ref{defn:S_2} that an arbitrary element of $S_2$ is a product of terms of the form $q(a)$ and terms of the form $p_2(\mu^i\eta_2\mu^{-i})$ for some  $a\in A$, $q(t)$ in $\Q[t,t^{-1}]$ where $\widetilde{(p_2,q)}= 1$, $q(1)\neq 0$, and $\mu$ generates $\pi/\pi^{(1)}$. Since $A$ is a torsion-free abelian group, \eqref{eq:storsion} may be viewed as an equation in $\Q F$ for some free abelian group $F\subset A$ of finite rank $r$. Since $\Q F$ is a UFD and since $\widetilde{(p_2,q)}= 1$ we can apply the following.

\begin{prop}[{\cite[Prop.4.5]{CHL5}}]\label{prop:characterizationequivalent} Suppose $p(t), q(t)\in \Q [t,t^{-1}]$ are non-zero. Then $p$ and $q$ are strongly coprime if and only if, for any finitely-generated free abelian group $F$ and any nontrivial $a,b\in F$, $p(a)$ is relatively prime to $q(b)$ in $\Q F$ (a unique factorization domain).
\end{prop}

\noindent Thus the greatest common divisor, in $\Q F$, of $p_2(\eta_1)$ and $q(a)$ is a unit (note that if $a$ is trivial in $F$ then $q(a)=q(1)\neq 0$ is itself a unit). Thus $p_2(\eta_1)$ divides the product of the terms of the form $p_2(\mu^i\eta_2\mu^{-i})$. Choose a basis, $\{x,x_2,\dots,x_r\}$, for $F$ in which $\eta_1=x^r$ for some $r>0$ (since $\eta_1\neq 0$ by Lemma~\ref{lem:easystep}) and $\mu^i\eta_2\mu^{-i}=x^{n_i}x_2^{n_{i,2}}\cdots x_r^{n_{i,r}}$. Then we may regard $\Q F$ as a Laurent polynomial ring in the variables $\{x,x_2,\dots,x_r\}$. Since $p_2$ is not zero and not a unit, there exists  a non-zero complex root $x=\tau$ of $p_2(x^r)$. Suppose $\tilde{p}(x)$ is an irreducible factor (in $\Q F$) of $p_2(x^r)$ of which $\tau$ is a root. Then, for some $i$, $\tilde{p}(x)$ divides $p_2(x^{n_i}x_2^{n_{i,2}}\cdots x_r^{n_{i,r}})$. Then $\tau$ must be a zero of $p_2(x^{n_i}x_2^{n_{i,2}}\cdots x_r^{n_{i,r}})$ for \emph{every} complex value of $x_2,\dots,x_r$. This is impossible unless each $n_{i,j}=0$. Thus, for this value of $i$, $\mu^i\eta_2\mu^{-i}=x^{n}$, in $F$, for some $n$. Note $n\neq 0$ since $\eta_2$ is nontrivial by Lemma~\ref{lem:easystep}. Thus
\begin{equation}\label{eq:in alex}
\mu^i\eta_2^r\mu^{-i}=(\mu^i\eta_2\mu^{-i})^r=x^{nr}=\eta_1^n,
\end{equation}
for some $i$ and some non-zero integers $n$ and $r$. This equation holds in $A$. However, the circles $\mu$, $\eta_2$ and $\eta_1$ all live in $M_{\mathfrak{R}^{m}}$ and in fact can be interpreted in $\mathcal{A}^{\Z}(\mathbb{E}^m)$ (the left-hand copy of $\mathbb{E}^m$). But recall that in the proof of Lemma~\ref{lem:easystep} we showed that the map
$$
\mathcal{A}^{\Z}(\mathbb{E}^m)\to \mathcal{A}^{\Z}(W)\to \frac{\pi^{(1)}}{\pi^{(2)}_\mathcal{S}}\equiv A
$$
is injective. Hence if \eqref{eq:in alex} holds in $A$ then it holds as an equation in $\mathcal{A}^{\Z}(\mathbb{E}^m)$, and hence also in $\mathcal{A}(\mathbb{E}^m)$, where, in module notation, it has the form
$$
(t_*)^i(r\eta_2)=n\eta_1.
$$
But the simple computation in the following Lemma proves that this is impossible.

\begin{lem}\label{lem:t-star} Let $m$ be a non-zero integer, let $\mathbb{E}^m$ be the knot of Figure~\ref{fig:Em} and let $\langle \eta_i \rangle, ~i=1,2$ be the subspace of $\mathcal{A}(\mathbb{E}^m)$ generated by the circle $\eta_i$ shown in Figure~\ref{fig:frakR}. Then, under the automorphism
$$
t_*:\mathcal{A}(\mathbb{E}^m)\to\mathcal{A}(\mathbb{E}^m),
$$
for every integer $k$, $(t_*)^k(\langle \eta_2 \rangle)\cap \langle \eta_1 \rangle= \vec{0}$.
\end{lem}
\begin{proof}We may assume that $m>0$. If $V$ is the Seifert matrix for $\mathbb{E}^m$ as in the proof of Proposition~\ref{prop:alexpolyscoprime}, with respect to the basis $\{a_i\}$ consisting of the cores of the obvious bands where $\ell k(a_i,\eta_i)=1$, then the rational Alexander module is presented by $V-tV^T$ with respect to the basis $\{\eta_1,\eta_2\}$ where the relations are given by the columns, that is, $(V-tV^T)\vec{v}=\vec{0}$ for all $\vec{v}$. Since $V$ has non-zero determinant, upon left multiplying the latter equation by $V^{-1}$, one recovers the fact that the automorphism $t^*$ is given by left multiplication by $(V^{-1})^TV$. Hence
$$
t^*=\frac{1}{m^2}\left(\begin{matrix} m^2+1 & m\cr m & m^2
\end{matrix}\right)=\frac{1}{m^2} M,
$$
for $M$ as indicated, with respect to the basis $\{\eta_1,\eta_2\}$. It then suffices to prove that, for any $k$, there is no non-zero solution $(x_0,y_0)$ to the equation
$$
M^k\left(\begin{matrix} 0 \cr y_0
\end{matrix}\right)=\left(\begin{matrix} x_0 \cr 0
\end{matrix}\right).
$$
If there were such a solution $(x_0,y_0)$ then there would be one with $y_0>0$. Let $B=\{(x,y)~|~x\geq 0,~y>0\}$. Since
$$
\left(\begin{matrix} m^2+1 & m\cr m & m^2
\end{matrix}\right)\left(\begin{matrix} x \cr y
\end{matrix}\right)=\left(\begin{matrix} (m^2+1)x+my \cr mx+m^2y
\end{matrix}\right)
$$
we observe that $M(B)\subset B$. But then if $k\geq 0$, $M^k(B)\subset B$. This is a contradiction since $(0,y_0)\in B$ but $(x_0,0)\notin B$.
Therefore there is no non-zero solution if $k\geq 0$. If $k<0$ then we have
$$
\left(\begin{matrix} 0 \cr y_0
\end{matrix}\right)=M^{-k}\left(\begin{matrix} x_0 \cr 0
\end{matrix}\right),
$$
where $-k=s>0$. As above if there were a non-zero solution then there would be one with $x_0>0$. Letting $A=\{(x,y)~|~x>0,~y\geq 0\}$, we observe that $M^{s}(A)\subset A$, $(x_0,0)\in A$ and $(0,y_0)\notin A$, which is a contradiction.
\end{proof}

This contradiction establishes ~\eqref{eq:nontrivialmodule}, finally finishing the proof of Lemma~\ref{lem:crucialstep}.
\end{proof}

We now complete the induction step in the proof of Proposition~\ref{prop:inductivestep}.

Suppose, for some $i$, $1\leq i\leq n-2$, Proposition~\ref{prop:inductivestep} holds, that is,  $\mu_{i}=\alpha_{i+1}$ is non-trivial, while $\overline{\mu}_{i}=\overline{\alpha}_{i+1}$ is trivial in
\begin{equation}\label{eq:indhypoth}
\frac{\pi_1(W_{i})^{(n-i)}}{\pi_1(W_{i})^{(n-i+1)}_{\mathcal{S}}}.
\end{equation}
To complete the inductive step we need to show that
\begin{equation}\label{eq:inducttriviality}
\overline{\mu}_{i-1}=\overline{\alpha}_{i}=0~\in
\frac{\pi_1(W_{i-1})^{(n-i+1)}}{\pi_1(W_{i-1})^{(n-i+2)}_{\mathcal{S}}}.
\end{equation}
and show that
\begin{equation}\label{eq:inductnontriviality}
\mu_{i-1}=\alpha_{i}\neq 0 ~\in
\frac{\pi_1(W_{i-1})^{(n-i+1)}}{\pi_1(W_{i-1})^{(n-i+2)}_{\mathcal{S}}}.
\end{equation}

\noindent By the inductive hypothesis and weak functoriality,
$$
\overline{\mu}_{i}\in \pi_1(W_{i})^{(n-i+1)}_{\mathcal{S}}\subset \pi_1(W_{i-1})^{(n-i+1)}_{\mathcal{S}}.
$$
But, by property ($1$) of Lemma~\ref{lem:mickeyfacts}, $\overline{\mu}_{i}\in \pi_1(M_{\overline{K}^i})$ normally generates $\pi_1(\overline{E}_i)$ so
$$
\pi_1(\overline{E}_i)\subset \pi_1(W_{i-1})^{(n-i+1)}_{\mathcal{S}},
$$
and so by property ($1$) of Proposition~\ref{prop:commseriesprops},
$$
[\pi_1(\overline{E}_i),\pi_1(\overline{E}_i)]\subset [\pi_1(W_{i-1})^{(n-i+1)}_{\mathcal{S}},\pi_1(W_{i-1})^{(n-i+1)}_{\mathcal{S}}]\subset\pi_1(W_{i-1})^{(n-i+2)}_{\mathcal{S}}.
$$
Since $\ell k(\overline{\alpha}_{i},\overline{R}^i)=0$,
$$
\overline{\alpha}_{i}\in [\pi_1(M_{\overline{K}^i}), \pi_1(M_{\overline{K}^i})]\subset [\pi_1(\overline{E}_i),\pi_1(\overline{E}_i)]\subset \pi_1(W_{i-1})^{(n-i+2)}_{\mathcal{S}}.
$$
This proves \eqref{eq:inducttriviality}.

Now to we need to prove \eqref{eq:inductnontriviality}. By property ($1$) of Lemma~\ref{lem:mickeyfacts}, the kernel of the map
$$
\pi_1(W_{i})\to \pi_1(W_{i}\cup E_{i}\cup \overline{E}_{i})= \pi_1(W_{i-1})
$$
is normally generated by the longitudes, $\ell_{i-1}, \overline{\ell}_{i-1}$, of the infecting knots $K^{i-1}$ and $\overline{K}^{i-1}$ viewed as curves in $S^3\setminus K^{i-1}\subset M_{K^{i}}\subset \partial W_{i}$ and $S^3\setminus \overline{K}^{i-1}\subset M_{\overline{K}^{i}}\subset \partial W_{i}$. But of course these lie in the second derived subgroups of $\pi_1(S^3\setminus K^{i-1})$ and $\pi_1(S^3\setminus \overline{K}^{i-1})$ respectively, and so lie in the second derived subgroups of $\pi_1(M_{K^{i}})$ and $\pi_1(M_{\overline{K}^{i}})$ respectively. But, by the induction hypothesis \eqref{eq:indhypoth},
\begin{equation}\label{eq:munminusone1}
\pi_1(M_{K^{i}})=\langle \mu_{i}\rangle\subset \pi_1(W_{i})^{(n-i)},
\end{equation}
and similarly for $\pi_1(M_{\overline{K}^{i}})$. It follows that both $\ell_{i-1}$ and $\overline{\ell}_{i-1}$ lie in
$$
\pi_1(W_{i})^{(n-i+2)}\subset \pi_1(W_{i})^{(n-i+2)}_{\mathcal{S}}.
$$
Thus the inclusion $W_{i}\to W_{i-1}$ induces an isomorphism
$$
\frac{\pi_1(W_{i})^{(n-i+1)}}{\pi_1(W_{i})^{(n-i+2)}_{\mathcal{S}}}\cong \frac{\pi_1(W_{i-1})^{(n-i+1)}}{\pi_1(W_{i-1})^{(n-i+2)}_{\mathcal{S}}},
$$
by weak functoriality and by \cite[Prop. 4.7]{CHL5}.

Consequently, to establish \eqref{eq:inductnontriviality}, it suffices to let $\pi=\pi_1(W_{i})$, and show that $\alpha_{i}$ is non-trivial in $\pi^{(n-i+1)}/\pi^{(n-i+2)}_{\mathcal{S}}$. Throughout the rest of the proof, we will abbreviate $W=W_{i}$, $\pi=\pi_1(W_{i})$, $J=K^{i}$ and $\overline{J}=\overline{K}^{i}$. Thus $M_J \subset \partial W$.

Consider the following commutative diagram (which we justify below) where $\G=\pi/\pi^{(n-i+1)}_\mathcal{S}$ and $\mathcal{R}=\Q\G S_{n-i+1}^{-1}$. Since $\alpha_{i}\in\pi_1(M_J)^{(1)}$ we have reduced \eqref{eq:inductnontriviality} to showing that $\alpha_{i}$ is \emph{not} in the kernel of the top row of the diagram.
$$
\begin{diagram}\label{diagram-pirelatehomology2}\dgARROWLENGTH=.8em
\node{\pi_1(M_J)^{(1)}}\arrow[2]{e,t}{j_*}\arrow{s,r}{\pi}\node[2]{\pi^{(n-i+1)}}\arrow{e,t}{\phi}\arrow{s,r}{} \node{\frac{\pi^{(n-i+1)}_\mathcal{S}}{\pi^{(n-i+2)}_\mathcal{S}}}\arrow{s,r}{j} \\
  \node{\mathcal{A}(J)\otimes \mathcal{R}}\arrow{e,t}{\cong }\node{ H_1(M_J;\mathcal{R})} \arrow{e,t}{j_*}\node{H_1(W;\R)}\arrow{e,t}
{\cong}\node{\frac{\pi^{(n-i+1)}_\mathcal{S}}{[\pi^{(n-i+1)}_\mathcal{S},\pi^{(n-i+1)}_\mathcal{S}]}\otimes \mathcal{R}}
\end{diagram}
$$
The $j_*$ in the upper row of the diagram is justified by \eqref{eq:munminusone1}. Now we consider the first map in the bottom row. By the inductive hypothesis \eqref{eq:inductnontriviality} the coefficient system $\pi\to \G$, when restricted to $\pi_1(M_J)$ is non-trivial:
$$
\pi_1(M_J)=\langle \mu_{i}\rangle \hookrightarrow \frac{\pi^{(n-i)}}{\pi^{(n-i+1)}_\mathcal{S}}\hookrightarrow \frac{\pi}{\pi^{(n-i+1)}_\mathcal{S}}\equiv \G,
$$
but also factors through $\pi_1(M_J)/\pi_1(M_J)^{(1)}\cong\Z$ because of \eqref{eq:munminusone1}. It follows that
$$
H_1(M_J;\Q\G)\cong H_1(M_J;\Q[t,t^{-1}])\otimes \Q\G\equiv \mathcal{A}(J)\otimes_{\Q[t,t^{-1}]} \Q\G,
$$
where $\Q[t,t^{-1}]$ acts on $\Q\G$ by $t\to \mu_{i}$. Hence
$$
H_1(M_J;\mathcal{R})\cong \mathcal{A}(J)\otimes \mathcal{R}.
$$
To justify the last map in the lower row, recall that $H_1(W;\Z\G)$ has an interpretation as the first homology module of the $\G$-covering space of $W$ corresponding to the kernel of $\pi\to \G$. Hence
$$
H_1(W;\Z\G)\cong \frac{\pi^{(n-i+1)}_\mathcal{S}}{[\pi^{(n-i+1)}_\mathcal{S},\pi^{(n-i+1)}_\mathcal{S}]}
$$
This completes the explanation of the diagram. Since, by Definitions~\ref{def:Sdefderivedlocalp} and \ref{defn:S_2},
$$
\pi^{(n-i+2)}_{\mathcal{S}}=\ker \left(\pi^{(n-i+1)}_\mathcal{S}\to \frac{\pi^{(n-i+1)}_\mathcal{S}}{[\pi^{(n-i+1)}_\mathcal{S},\pi^{(n-i+1)}_\mathcal{S}]}\to \frac{\pi^{(n-i+1)}_\mathcal{S}}{[\pi^{(n-i+1)}_\mathcal{S},\pi^{(n-i+1)}_\mathcal{S}]}\otimes \Q\G S_{n-i+1}^{-1}\right),
$$
it follows that the vertical map $j$ (in the diagram) is injective. Hence, to establish \eqref{eq:inductnontriviality}, it suffices to show that the class represented by $\alpha_{i}\otimes 1$ is \emph{not} in the kernel of the bottom row of the diagram.

Recall that $J\equiv K^{i}\equiv R^{i}_{\alpha_{i}}(K^{i-1})$  where $\alpha_{i}$ generates $\mathcal{A}(R^{i})$. Therefore $\mathcal{A}(J)\cong \mathcal{A}(R^i)$. By the hypotheses of Theorem~\ref{thm:main}, the Alexander polynomial of $R^{i}$ is $q_{i}(t)=p_{n-i+1}(t)$. Thus
\begin{equation}\label{eq:module}
\langle \alpha_{i}\otimes 1\rangle\cong\mathcal{A}(J)\otimes\mathcal{R}\cong \left (\frac{\Q\G}{p_{n-i+1}(\mu_i)\Q\G}\right )S_{p_{n-i+1}}^{-1}.
\end{equation}
where the last equality holds because, since $1\leq i\leq n-2$, it follows that $3\leq n-i+1\leq n$, so $S_{n-i+1}=S_{p_{n-i+1}}$, by Definition~\ref{defn:S_2}.

Suppose that $\alpha_{i}\otimes 1$ \emph{were} in the kernel of the bottom row of the diagram. We shall reach a contradiction. Recall that
$$
W=W_{i}\equiv V\cup E_n\cup E_{n-1}\cup \overline{E}_{n-1}\cup\dots\cup E_{i+1}\cup \overline{E}_{i+1}.
$$
Recall also that $V$ is an $(n.5,\mathcal{P})$-solution. Thus, by \eqref{eq:PinPprime2}, $V$ is an $(n.5,\mathcal{S})$-solution and, since $n-i+1\leq n$, $V$ is also an $(n-i+1,\mathcal{S})$-solution. One easily checks that
$$
\frac{H_2(W_{i})}{i_*(H_2(\partial W_{i}))}\cong H_2(V).
$$
Hence this group has a basis consisting of surfaces that satisfy parts ($2$) and ($3$) of Definition~\ref{def:Gnsolvable} (with $n-i+1$). Thus $W_i$ is an $(n-i+1,\mathcal{S})$-bordism (\cite[Definition 7.11]{CHL5}). By ~\cite[Thms. 7.14, 7.15]{CHL5}, if $P$ is the kernel of the map
$$
j_*:H_1(M_J;\mathcal{R})\to H_1(W;\mathcal{R}),
 $$
then $P$ is isotropic for the Blanchfield linking form on $H_1(M_J;\mathcal{R})$. Therefore if the generator $\alpha_{i}\otimes 1$ were in $P$, it would follow that this Blanchfield form were identically zero on $H_1(M_J;\mathcal{R})$. But by \cite[Lemma 7.16]{CHL5} this form is non-singular. This would imply that $H_1(M_J;\mathcal{R})=0$. This is a contradiction once we show that \eqref{eq:module} is in fact a non-trivial module. It is shown in \cite[Theorem 4.12]{CHL5} that
\begin{equation}
\frac{\Q\G}{p_{n-i+1}(\mu_i)\Q\G}\hookrightarrow \left (\frac{\Q\G}{p_{n-i+1}(\mu_i)\Q\G}\right )S_{p_{n-i+1}}^{-1},
\end{equation}
is a monomorphism (using that $p_{n-i+1}(t)\neq 0$ and that $\mu_i$ lies in the abelian normal subgroup $A=\pi^{(n-i)}/\pi^{(n-i+1)}_\mathcal{S}\subset \G$). This reduces us to showing that
\begin{equation}\label{eq:module2}
\frac{\Q\G}{p_{n-i+1}(\mu_i)\Q\G}\neq 0.
\end{equation}
By the hypotheses of Theorem~\ref{thm:main}, $p_{n-i+1}(t)=q_{i}(t)$ is not a unit. The map $\Z\to \G$ given by $t\mapsto \mu_i$ is not zero by the inductive hypothesis \eqref{eq:indhypoth}. Thus $\langle \mu_i\rangle\subset\G$ and $\Q\G$ is a free $\Q[\mu_i,\mu_i^{-1}]$-module on the cosets of $\langle\mu_i\rangle\in \G$. In the same manner as we showed earlier in the proof, it follows that $p_{n-i+1}(\mu_i)$ is not a unit in the domain $\Q\G$. Therefore \eqref{eq:module2} holds.

This finishes, finally, the inductive step and hence the entire proof of Proposition~\ref{prop:inductivestep}, which in turn completes the proofs of \eqref{eq:nontrivialonK} and \eqref{eq:trivialonKbar}.
\end{proof}
\end{subsection}

Having established \eqref{eq:PinPprime}, \eqref{eq:nontrivialonK} and \eqref{eq:trivialonKbar}, the proof of Theorem~\ref{thm:main} is complete.

\end{proof}

\bibliographystyle{plain}
\bibliography{mybib6mathscinet,mybib7mathscinet}
\end{document}